\documentclass{fundam}        

\usepackage{chemarr}
\usepackage{amsfonts}
\usepackage{amssymb}
\usepackage{latexsym}
\usepackage{graphicx}
\usepackage{amsbsy}
\usepackage{psfrag}
\usepackage{color}

\newcommand{\D}{\displaystyle}
\newcommand{\restr}[1]{\raisebox{-0.3em}{$\lb|_{#1}\rb.$}} 
\newcommand{\ignore}[1]{}   
\newcommand{\Caption}[2]{{\footnotesize\caption{\label{#1} #2}}}   
\newcommand{\breath}{\medskip} 
\newtheorem{thm}{Theorem}[section]

\newcounter{claimcount}[thm]  
\newcounter{subclaimcount}[claimcount]
\newtheorem{prop}[thm]{Proposition}   

\newtheorem{cor}[thm]{Corollary}

\newcommand{\dfn}{\sf\em} 
\newcommand{\Theorem}[2]{\begin{thm}{\sf #1}  #2 \end{thm}}
\newcommand{\Proposition}[2]{\begin{prop}{\sf #1}  #2 \end{prop}}
\newcommand{\Lemma}[2]{\begin{lemma}{\sf #1}  #2 \end{lemma}}
\newcommand{\Corollary}[2]{\begin{cor}{\sf #1}  #2 \end{cor}} 
\newcommand{\thmfont}[1]{{\sl #1}}     
\newcommand{\EXAMPLE}[1]{        \refstepcounter{thm}                     \begin{list}{} 			{\setlength{\leftmargin}{0em} 			\setlength{\rightmargin}{0em}}        \item {\bf Example \thethm:} #1                   \hfill$\diamondsuit$  \end{list}   			}     
\newcommand{\bthmlist}{ \begin{list}{{\bf(\alph{enumi})}} {\usecounter{enumi} \setlength{\leftmargin}{0.5em} \setlength{\itemsep}{0.2em} \setlength{\topsep}{0.2em} \setlength{\itemindent}{0em} \setlength{\parsep}{0em} \setlength{\rightmargin}{0em}} } 
\newcommand{\ethmlist}{\end{list}}    
\newcommand{\Claim}[1]{\refstepcounter{claimcount}                \noindent {\sc Claim \theclaimcount: \ }\thmfont{ #1}} 
\newcommand{\Subclaim}[1]{\refstepcounter{subclaimcount}                \noindent {\sc Claim \theclaimcount.\thesubclaimcount: \ }\thmfont{ #1}} 
\newcommand{\subclaim}{\Subclaim}
\newcommand{\bprf}[1][Proof:]{\begin{list}{} 			{\setlength{\leftmargin}{0.7em} 			\setlength{\rightmargin}{0em} 			\setlength{\listparindent}{1em}}                         \item {\em \hspace{-1.1em}  #1  }} 
\newcommand{\eprf}{\end{list}} 
\newcommand{\bthmprf}{\bprf}
\newcommand{\bclaimprf}{\bprf}
\newcommand{\bsubclaimprf}{\bprf}
\newcommand{\ethmprf}{ \hfill$\Box$  \eprf  \breath  } 
\newcommand{\eclaimprf}{ \hfill $\Diamond$~{\scriptsize {\tt Claim~\theclaimcount}}\eprf}  
\newcommand{\esubclaimprf}{ \hfill $\triangledown$~{\scriptsize  {\tt Claim~\theclaimcount.\thesubclaimcount}}\eprf}  
\newcommand{\qed}{\QED}     
\newcommand{\beq}{\begin{eqnarray*}}
\newcommand{\eeq}{\end{eqnarray*}} 
\newcommand{\beqn}{ \begin{equation} }
\newcommand{\eeqn}{ \end{equation} }
\newcommand{\blist}{\begin{enumerate}}
\newcommand{\elist}{\end{enumerate}} 
\newcommand{\bitem}{\begin{itemize}}
\newcommand{\eitem}{\end{itemize}} 
\newcommand{\bdesc}{\begin{description}}
\newcommand{\edesc}{\end{description}}   
\newcommand{\If}{\mbox{\ if \ }}  
\newcommand{\Folner}{F{\o}lner }
\newcommand{\Cesaro}{Ces\`aro }
\newcommand{\done}{{\mathsf{ 1\!\!1}}} 
\newcommand{\dB}{{\mathbb{B}}}
\newcommand{\dC}{{\mathbb{C}}}
\newcommand{\dD}{{\mathbb{D}}}
\newcommand{\dE}{{\mathbb{E}}}
\newcommand{\dF}{{\mathbb{F}}}
\newcommand{\dG}{{\mathbb{G}}}
\newcommand{\dH}{{\mathbb{H}}}
\newcommand{\dK}{{\mathbb{K}}}
\newcommand{\dN}{{\mathbb{N}}}
\newcommand{\dP}{{\mathbb{P}}}
\newcommand{\dQ}{{\mathbb{Q}}}
\newcommand{\dR}{{\mathbb{R}}}
\newcommand{\dT}{{\mathbb{T}}}
\newcommand{\dU}{{\mathbb{U}}}
\newcommand{\dV}{{\mathbb{V}}}
\newcommand{\dW}{{\mathbb{W}}}
\newcommand{\dX}{{\mathbb{X}}}
\newcommand{\dY}{{\mathbb{Y}}}
\newcommand{\dZ}{{\mathbb{Z}}}       %  Define Bar macros 
\newcommand{\barE}{{\overline{E}}}
\newcommand{\bL}{{\mathbf{ L}}}
\newcommand{\bT}{{\mathbf{ T}}}
\newcommand{\ba}{{\mathbf{ a}}}
\newcommand{\bb}{{\mathbf{ b}}}
\newcommand{\bc}{{\mathbf{ c}}}
\newcommand{\bd}{{\mathbf{ d}}}
\newcommand{\be}{{\mathbf{ e}}}
\newcommand{\bg}{{\mathbf{ g}}}
\newcommand{\bi}{{\mathbf{ i}}}
\newcommand{\bs}{{\mathbf{ s}}}
\newcommand{\bt}{{\mathbf{ t}}}
\newcommand{\bw}{{\mathbf{ w}}}
\newcommand{\bx}{{\mathbf{ x}}}
\newcommand{\by}{{\mathbf{ y}}}
\newcommand{\bdel }{{\boldsymbol{\delta}}}
\newcommand{\blam }{{\boldsymbol{\lambda}}}
\newcommand{\sA}{{\mathcal{ A}}}
\newcommand{\sB}{{\mathcal{ B}}}
\newcommand{\sC}{{\mathcal{ C}}}
\newcommand{\sD}{{\mathcal{ D}}}
\newcommand{\sE}{{\mathcal{ E}}}
\newcommand{\sF}{{\mathcal{ F}}}
\newcommand{\sK}{{\mathcal{ K}}}
\newcommand{\sL}{{\mathcal{ L}}}
\newcommand{\sO}{{\mathcal{ O}}}
\newcommand{\sP}{{\mathcal{ P}}}
\newcommand{\sX}{{\mathcal{ X}}}

\newcommand{\gA}{{\mathfrak{ A}}}
\newcommand{\gB}{{\mathfrak{ B}}}
\newcommand{\gC}{{\mathfrak{ C}}}
\newcommand{\gD}{{\mathfrak{ D}}}
\newcommand{\gE}{{\mathfrak{ E}}}

\newcommand{\gG}{{\mathfrak{ G}}}
\newcommand{\gM}{{\mathfrak{ M}}}
\newcommand{\gO}{{\mathfrak{ O}}}
\newcommand{\gS}{{\mathfrak{ S}}}
\newcommand{\gT}{{\mathfrak{ T}}}
\newcommand{\gU}{{\mathfrak{ U}}}
\newcommand{\gW}{{\mathfrak{ W}}}
\newcommand{\gX}{{\mathfrak{ X}}}
\newcommand{\gh}{{\mathfrak{ h}}}
\newcommand{\go}{{\mathfrak{ o}}}
\newcommand{\alp }{\alpha}
\newcommand{\bet }{\beta}
\newcommand{\gam }{\gamma}
\newcommand{\del }{\delta}
\newcommand{\eps }{\epsilon}
\newcommand{\kap }{\kappa}
\newcommand{\lam }{\lambda}
\newcommand{\sig }{\sigma} 
\newcommand{\omg }{\omega}
\newcommand{\Del }{\Delta}
\newcommand{\h}{\widehat} 
\newcommand{\fZ}{{\mathsf{ Z}}} 
\newcommand{\fe}{{\mathsf{ e}}}
\newcommand{\ff}{{\mathsf{ f}}}
\newcommand{\fh}{{\mathsf h}} 
\newcommand{\fp}{{\mathsf{ p}}}
\newcommand{\fu}{{\mathsf{ u}}}
\newcommand{\fv}{{\mathsf{ v}}}
\newcommand{\fx}{{\mathsf{ x}}}
\newcommand{\fy}{{\mathsf{ y}}}
\newcommand{\fz}{{\mathsf{ z}}}    %  Define Tilde macros 
\newcommand{\tl}{\widetilde} 
\newcommand{\undba}{{\underline{\mathbf{ a}}}}
\newcommand{\lb}{\left}
\newcommand{\rb}{\right} 
\newcommand{\maketall}{\rule[-0.5em]{0em}{1em}}        
\newcommand{\map}{{\longrightarrow}}
\newcommand{\goto}{{\rightarrow}}
\newcommand{\into}{{\map}}
\newcommand{\statement}[1]{\lb(  \maketall       \begin{minipage}{40em}       \begin{tabbing}         #1        \end{tabbing}      \end{minipage}  \rb)}     
\newcommand{\oo}{{\infty}}        
\newcommand{\X}{\times}
\newcommand{\x}{\X}
\newcommand{\dirsum}{\oplus}
\newcommand{\Dirsum}{\bigoplus} 
\newcommand{\union}{\cup}
\newcommand{\Union}{\bigcup}
\newcommand{\intsct}{\cap}
\newcommand{\Intsct}{\bigcap}
\newcommand{\disj}{\sqcup}
\newcommand{\Disj}{\bigsqcup}   
\newcommand{\set}[2]{{\left\{ #1 \; ; \; #2 \right\} }} 
\newcommand{\Id}[1]{{\mathbf{ Id}_{{#1}}}}
\newcommand{\chr}[1]{{{\done}_{{#1}}}} 
\newcommand{\choice}[1]{{\lb\{ \begin{array}{rcl}                                 #1                                \end{array}  \rb.  }}                     %        
\newcommand{\eeequals}[1]{\raisebox{-0.9ex}{$\overline{\overline{{\scriptscriptstyle{\mathrm{#1}}}}}$}} 
\newcommand{\suuupseteq}[1]{\raisebox{-1.3ex}{$\stackrel{\D\supseteq}{\scriptscriptstyle{\mathrm{#1}}}$}} 
\newcommand{\closeto}[1]{{{\raisebox{-1ex}    {$\widetilde{\ {\scriptstyle #1}\ }$}}}}      
\newcommand{\Fix}[1]{{\sf Fix}\lb[#1\rb]}
\newcommand{\shift}[1]{\sig^{#1}}    
\newcommand{\diam}[1]{{\sf diam}\lb[#1\rb]} 
\newcommand{\Aut}[2][]{\mathsf{Aut}_{#1} \lb(#2\rb)}
\newcommand{\Real}{\dR}
\newcommand{\Natur}{\dN}
\newcommand{\Zahl}{\dZ}
\newcommand{\Zahlmod}[1]{{\Zahl_{/#1}}}
\newcommand{\Rat}{\dQ}
\newcommand{\Cplx}{\dC}
\newcommand{\Torus}[1]{{{\dT}^{#1}}} 
\newcommand{\CC}[1]{{\lb[ #1 \rb]}}
\newcommand{\CO}[1]{{\lb[ #1 \rb)}}
\newcommand{\OC}[1]{{\lb( #1 \rb]}}
\newcommand{\IMPLIES}{\ensuremath{\Longrightarrow}}
\newcommand{\AND}{\mbox{\ and \ }}   
\newcommand{\ZD}[1][D]{{\Zahl^{#1}}}
\newcommand{\RD}[1][D]{{\Real^{#1}}}
\newcommand{\AZD}[1][D]{\sA^{\ZD[#1]}}
\newcommand{\AZ}{\sA^{\Zahl}}
\newcommand{\BZD}[1][D]{\sB^{\ZD[#1]}}
\newcommand{\tlfz}{\widetilde{\fz}}
\newcommand{\tlfZ}{\widetilde{\fZ}}
\newcommand{\tlfe}{\widetilde{\fe}}       
\newcommand{\bone}{\boldsymbol{1}}     
\newcommand{\black}{{\scriptstyle{\blacksquare}}}
\newcommand{\white}{{\scriptstyle{\square}}}  
\newcommand{\Sft}{\gA}
\newcommand{\SFT}{(\Sft,\Phi,\shift{})} 
\newcommand{\Cae}{\sC_{\ae}} 
\newcommand{\Spec}[1][\Sft,\Phi,\shift{}]{\mathrm{S}^{\scriptscriptstyle\mathrm{pec}}\!\lb(#1\rb)}
\newcommand{\Specae}[1][\Sft,\Phi,\shift{}]{\mathrm{S}^{\scriptscriptstyle\mathrm{pec}}_{\scriptscriptstyle\ae}\!\lb(#1\rb)}
\newcommand{\hSpecae}[1][\Sft,\Phi,\shift{}]{\widehat{\mathrm{S}^{\scriptscriptstyle\mathrm{pec}}_{\scriptscriptstyle\ae}}\!\lb(#1\rb)}
\newcommand{\RatSpec}[1][\Sft,\Phi,\shift{}]{\mathrm{S}^{\scriptscriptstyle\mathrm{pec}}_{\scriptscriptstyle\Rat}\!\lb(#1\rb)}
\newcommand{\hRatSpec}[1][\Sft,\Phi,\shift{}]{\widehat{\mathrm{S}^{\scriptscriptstyle\mathrm{pec}}_{\scriptscriptstyle\Rat}}\!\lb(#1\rb)}  
\newcommand{\Eigenspace}[1][\blam]{\sE^{_{\!i\!g\!e\!n}}_{#1}}  
\newcommand{\NRemark}[1]{                     \begin{list}{} 			{\setlength{\leftmargin}{0em} 			\setlength{\rightmargin}{0em}}        \item        \refstepcounter{thm} {\bf Remark \thethm:} #1                    \hfill$\diamondsuit$  \end{list}   			} 
\newcommand{\Remarks}[1]{                      \begin{list}{} 			{\setlength{\leftmargin}{0em} 			\setlength{\rightmargin}{0em}}         \item {\bf Remarks:} #1                   \hfill$\diamondsuit$  \end{list}   			} 
\newcommand{\Remark}[1]{                      \begin{list}{} 			{\setlength{\leftmargin}{0em} 			\setlength{\rightmargin}{0em}}         \item {\bf Remark:} #1                   \hfill$\diamondsuit$  \end{list}   			} 
\newcommand{\hgA}{\widehat{\gA}}
\newcommand{\tlgA}{\widetilde{\gA}}
\newcommand{\Nh}{\dH}
\newcommand{\nh}{\fh}
\newcommand{\Kurka}{K\r{u}rka } 
\newcommand{\Mono}{\gM{\scriptstyle\go}}
\newcommand{\Checker}{\gC{\scriptstyle\gh}}
\newcommand{\Dom}{\mathfrak{D{\scriptstyle\!om}}}
\newcommand{\adjacent}[1][]{\stackrel{#1}{\leadsto}}   
\newcommand{\ECA}[1]{{}_{\varepsilon}\!\Phi_{\mbox{\tiny #1}}}
\newcommand{\energy}[1][\ba]{\sF_{#1}} 
\newcommand{\unflawed}[1][r]{\dG_{#1}}
\newcommand{\displace}[1]{\del\raisebox{-0.3em}{\tiny $\stackrel{-\!\!\!-\!\!\!-\!\!\!\rightharpoonup}{#1}$}}     
\newcommand{\parag}[1]{\vspace{0.5em}   \paragraph*{{\em #1}}}   
\newcommand{\Adm}[1][r]{\gA_{(#1)}}  
\newcommand{\projdisl}[1]{} 
\newcommand{\Good}[1][]{\hgA_{#1}} 
\begin{document}

%\ignore
{
\title{Spectral domain boundaries in cellular automata} 

\author{Marcus Pivato\thanks{Research partially supported by NSERC Canada, and by the Van Vleck Fund of Wesleyan University.}
\\
 Dept. of Mathematics \& Computer Science, Wesleyan University,
 Middletown CT, USA \\
 and Dept. of Mathematics, Trent University, 
Peterborough, Ontario,  Canada \\
{\tt pivato@xaravve.trentu.ca}
}
\address{Dept. of Mathematics, Trent University, 2151 East Bank Drive,
Peterborough, Ontario K9L 1Z6,  Canada}

\runninghead{M. Pivato}{Spectral domain boundaries in cellular automata} 

\newcounter{volume}
\setcounter{volume}{77}
\issue{\Roman{volume} (2007)}

\maketitle

\begin{abstract}
Let $\AZD$ be the Cantor space of $\ZD$-indexed configurations in a
finite alphabet $\sA$, and let $\shift{}$ be the $\ZD$-action of
shifts on $\AZD$.  A {\dfn cellular automaton} is a continuous,
$\shift{}$-commuting self-map $\Phi$ of $\AZD$, and a {\dfn $\Phi$-invariant
subshift} is a closed, $(\Phi,\shift{})$-invariant subset
$\gA\subset\AZD$.  Suppose $\ba\in\AZD$ is $\gA$-admissible everywhere
except for some small region we call a {\dfn defect}. It has been
empirically observed that such defects persist under iteration of
$\Phi$, and often propagate like `particles' which coalesce or
annihilate on contact.  We use spectral theory to
explain the persistence of some defects
under $\Phi$, and partly explain the outcomes of their collisions.

{\footnotesize
\breath
\ignore{
\begin{tabular}{rl}
{\bf MSC:}& 37B15 (primary), 68Q80 (secondary)\\
{\bf Keywords:}& Cellular automaton, subshift, defect, kink, domain boundary.
\end{tabular}}
}
\end{abstract}
  An often-observed phenomenon in cellular automata is the emergence
and persistence of homogeneous `domains' (each characterized by a
particular spatial pattern), separated by {\em defects} (analogous to
`domain boundaries' or `kinks' in a crystalline solid) which evolve
over time, propagating and occasionally colliding.  Such defects were
first empirically observed by Grassberger in the `elementary' cellular
automata  or `ECA' (radius-one CA on $\{0,1\}^\Zahl$) with numbers \#18,
\#122, \#126, \#146, and \#182 \cite{Gra84,Gra83} and also noted in ECA \#184,
considered as a simple model of surface growth \cite[\S III.B]{KrSp88}.
 Later, Boccara {\em et al.} empirically investigated 
the motion and interactions of defects in ECA
\#18, \#54, \#62, and \#184 (see Figure \ref{fig:defect.intro}),
and longer range totalistic CA \cite{BNR91,BoRo91}; see also 
\cite[\S3.1.2.2 \& \S3.1.4.4]{Ilachinski}.
Eloranta made the first rigorous mathematical study of defects in
\cite{Elo93a,Elo93b,Elo94,Elo95,ElNu}, while  Crutchfield and Hanson 
developed an empirical methodology 
called {\em Computational Mechanics} 
\cite{CrHa92,CrHa93a,CrHa93b,CrHa97,CrHR,Han}.
In another paper \cite{PivatoDefect0}, we characterized
the propagation of defects in one-dimensional subshifts of finite type.

A mathematical theory of cellular automaton defect dynamics is starting
to emerge, but some basic questions remain unanswered:
\blist
  \item \label{Q:defn}
 What is the right definition of `defect?' What constitutes a `regular domain'?

  \item \label{Q:persist} Why should a defect persist under the action
of cellular automata, rather than disappearing?  Are there
`topological' constraints imposed by the structure of the underlying
domain, which make defects indestructible?

  \item \label{Q:chem} When defects collide, they often coalesce into
a new type of defect, or mutually annihilate.  Is there a `chemistry'
governing these defect collisions?

  \item \label{Q:inv} Can we assign algebraic invariants to defects,
which reflect {\bf(a)} the `topological constraints' of question
\#\ref{Q:persist} or {\bf(b)} the `defect chemistry' of question
\#\ref{Q:chem}?
\elist
  This paper is organized as follows: in \S\ref{S:defect} we propose
an answer to question \#\ref{Q:defn}, which is a synthesis of
several previous approaches (see below).
We also introduce several examples which recur throughout the
paper.  In \S\ref{S:interface} and \S\ref{S:disl}, we use spectral
theory to address question \#\ref{Q:inv} for two types of domain
boundaries: {\em interfaces} and {\em dislocations}.  We develop
spectral invariants which answer question \#\ref{Q:persist},
and partially answer question \#\ref{Q:chem}.

\begin{figure}[h]
\centerline{\begin{tabular}{ll}
\includegraphics[clip=true, trim=0 0 300 240, width=17.5em,height=6.5em]{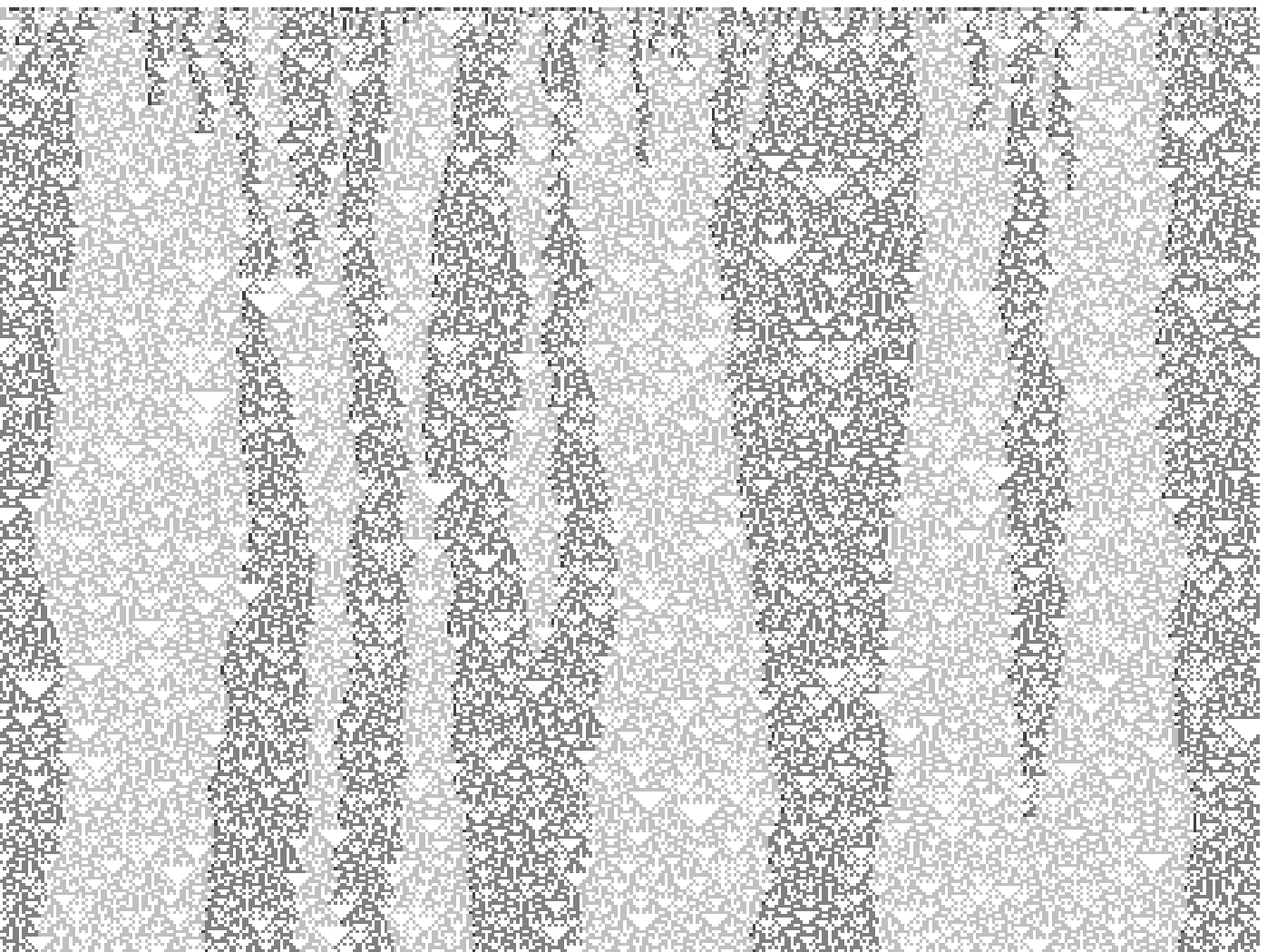} 
&\includegraphics[clip=true, trim=0 0 350 30, width=17.5em,height=6.5em]{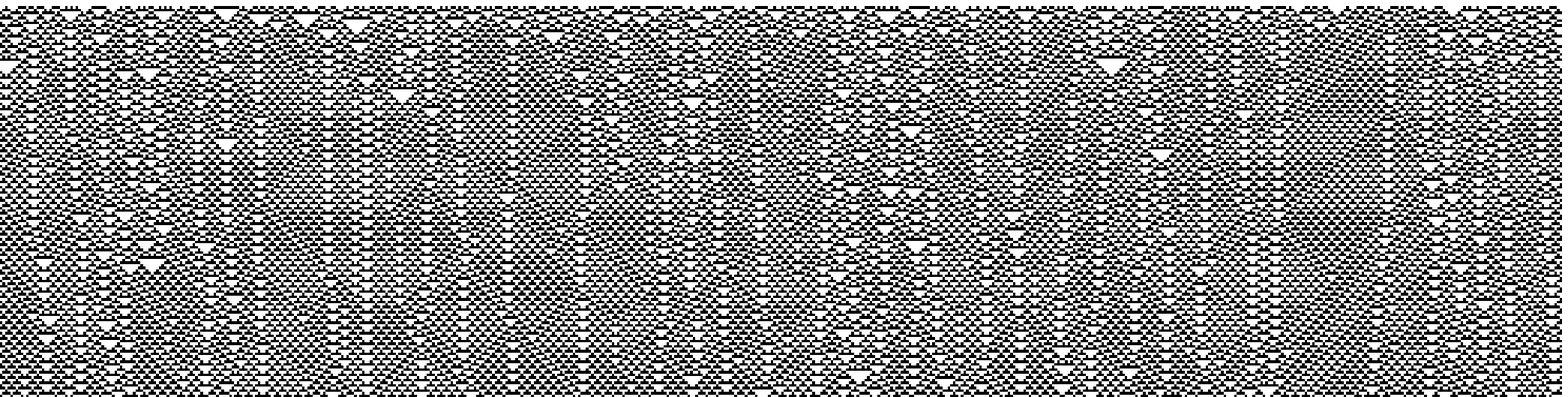} \\
{\footnotesize{\bf(A)} \ ECA \#18 (colourized to show defects)}& {\footnotesize{\bf(B)} \ ECA \#54 }\\
\includegraphics[clip=true, trim=0 30 350 0, width=17.5em,height=6.5em]{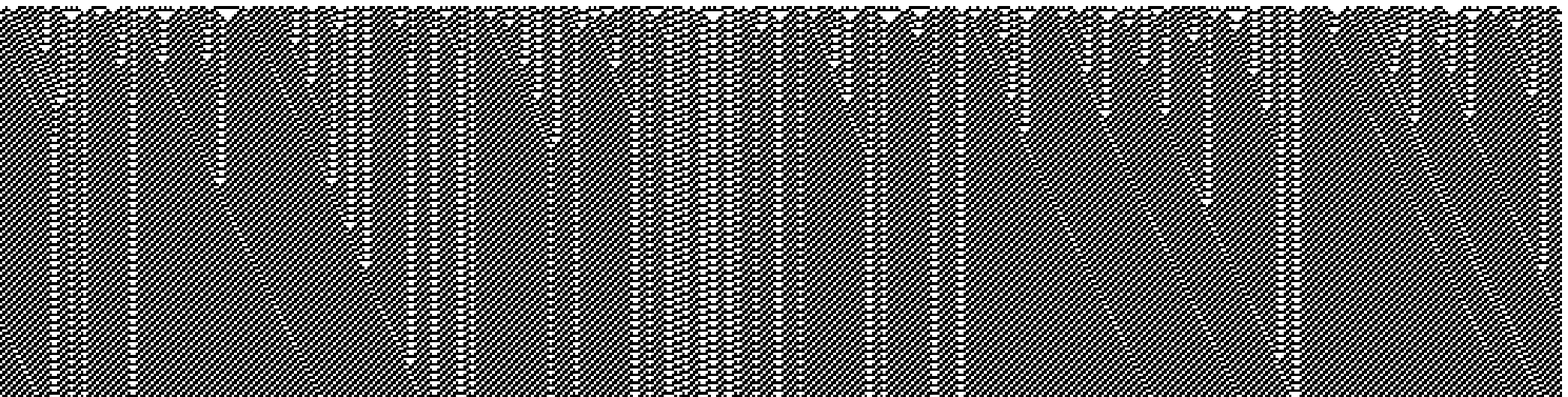} 
& \includegraphics[clip=true, trim=0 30 350 0, width=17.5em,height=6.5em]{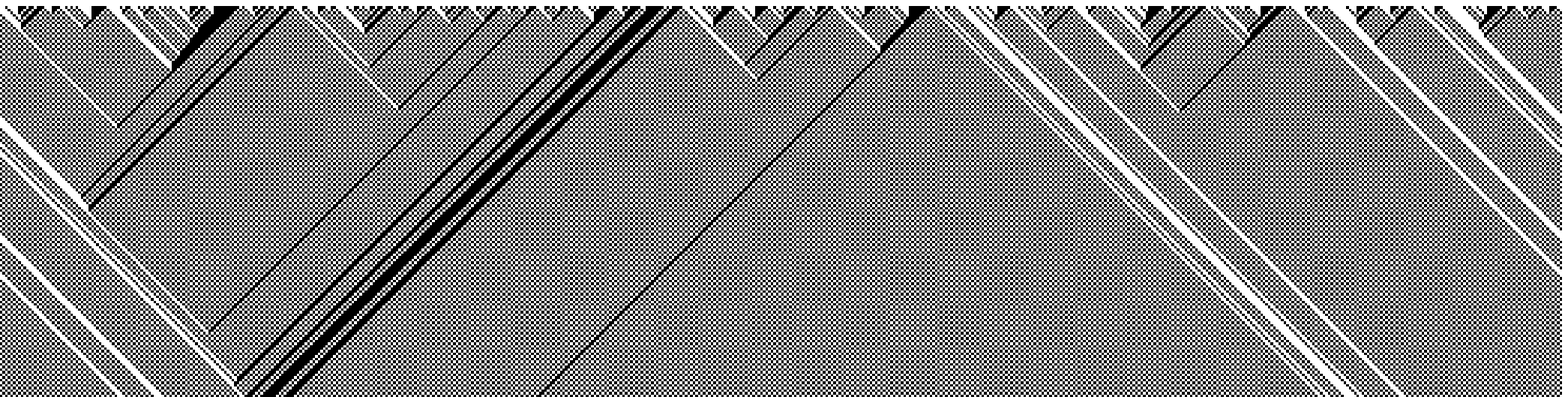} \\
{\footnotesize{\bf(C)} \ ECA \#62 } & {\footnotesize{\bf(D)} \ ECA \#184} \\
\end{tabular}}
\Caption{fig:defect.intro}{
 Spacetime diagrams showing defect dynamics in one-dimensional
cellular automata.  Each picture show 120 timesteps on a 250 pixel array
(time increases downwards).}
\end{figure}

\parag{Background:}
Roughly speaking, there are three  approaches to question \#\ref{Q:defn}:
\bitem
  \item[{\bf(a)}] `Domains' are identified with `invariant subalphabets' of the CA.
`Defects' are transitions from one subalphabet to another.

  \item[{\bf(b)}] `Defects' are identified by evaluating some numerical `weight
function' on the subwords of a sequence.
`Domains' are  regions where this weight function takes the value
0. 

  \item[{\bf(c)}] `Domains' are identified with some subshift (finite type
or sofic).  `Defects' are the (minimal-length) `forbidden words' of
this subshift.
\eitem

Approach {\bf(a)} was first developed in
\cite{ElNu} to prove Lind's conjecture \cite[\S5]{Lin84} that
the defects of ECA\#18 perform random walks.  It was then applied
to ECA's \#22, \#54, \#184 and other one-dimensional CA
\cite{Elo93a,Elo93b}, and extended  to defect ensembles \cite{Elo94}
and two-dimensional domain boundaries \cite{Elo95}.  Approach {\bf(b)}
was proposed by \Kurka \cite{Kur03}, and used to verify another 
conjecture of Lind \cite[\S5]{Lin84}, that some
cellular automata (such as ECA\#18) converge `in measure' to certain
limit sets through a process of defect coalescence/annihilation; see
also \cite{Kur05,KuMa00,KuMa02} for related ideas.  Approach {\bf(c)}
was first suggested in \cite[\S5]{Lin84}, and later elaborated in
\cite{CrHa92,CrHa93a,CrHa93b,CrHa97,CrHR,Han}, where  each `domain' 
was identified with a regular language (or equivalently, with a sofic
shift), which could be digitally filtered out of the spacetime diagram
using a finite automaton, thereby revealing the defect trajectories.
Approach {\bf(c)} was also used in \cite{PivatoDefect0}.

  Each approach has advantages and disadvantages.  Approach {\bf(b)}
is the most flexible, but also seems the most artificial, since we
must explicitly define the defects. In contrast, approach {\bf(c)}
allows defects to arise out of a simple and naturally occuring
background.  However, both {\bf(b)} and {\bf(c)} suffer from a problem
of defect `delocalization'.  A sofic shift cannot be characterized by
any finite set of forbidden words.  Thus, there must exist
huge yet indecomposable `defect particles'.  If
we arbitrarily assign these unwieldy defects a `location' at some point
in $\ZD$, then we are confronted with instantaneous
long-range interactions between defects, or the apparent
`teleportation' of defects through space.  To avoid this problem, we
could only admit subshifts of finite type as regular domains; then the
minimal-length forbidden words are bounded in length, so that defects
can be localized.  However, for certain CA (e.g. ECA\#18),
the regular domain is {\em not} of finite type ---it is strictly sofic
---hence these CA would escape our analysis.

  Approach {\bf(a)} avoids the `delocalization' problem.  For example,
the `defect-free' subshift of ECA\#18 is a sofic subshift [see Example
\ref{X:defect.energy}(c)], but when recoded using the alphabet
$\{00,01,10,11\}$, it is simply the disjoint union of two full shifts:
$\gE=\sE^\Zahl$ and $\gO=\sO^\Zahl$, where $\sE:=\{00,01\}$ and
$\sO:=\{00,10\}$ (the symbol $11$ is forbidden).  A defect is then a
transition from an $\sE$-valued sequence to a $\sO$-valued sequence.
However, {\bf(a)} is only appropriate for codimension-one defects
(i.e. `domain boundaries') and not for defects of codimension two
(e.g. `holes' in $\Zahl^2$, `strings' in $\Zahl^3$, etc.) or higher.
Even in codimension one, not every domain/defect problem can be
recoded in terms of invariant subalphabets, without obscuring
important information.

  In \S\ref{S:defect}, we will propose a definition which combines
features of {\bf(a)}, {\bf(b)}, and {\bf(c)}, and which is applicable
defects of any codimension, in any kind of subshift (sofic or otherwise).

\subsubsection*{Preliminaries \& Notation}

For any $L\leq R\in\Zahl$, we define 
\[
\CC{L...R} \ := \ \{L, L\!+\!1,\,\ldots,R\},\qquad
\CO{L...R} \ := \ \CC{L\ldots R\!-\!1}, \qquad 
\OC{L...R} \ := \ \CC{L\!+\!1\ldots R},\quad\mbox{etc.}
\] 
 We likewise define $\OC{-\oo....R}$, \ $\CO{L...\oo}$, etc.
Let $\sA$ be a finite alphabet.  Let $D\geq 1$, let $\ZD$ be the
$D$-dimensional lattice, and let $\AZD$ be the set of all
$\ZD$-indexed {\dfn configurations} of the form
$\ba=[a_{\fz}]_{\fz\in\ZD}$, where $a_\fz\in\sA$ for all $\fz\in\ZD$.
The {\dfn Cantor metric} on $\AZD$ is defined by
$d(\ba,\bb)=2^{-\Del(\ba,\bb)}$, where
$\Del(\ba,\bb):=\min\set{|\fz|}{a_\fz\neq b_\fz}$.  It follows that
$(\AZD,d)$ is a Cantor space (i.e. a compact, totally disconnected,
perfect metric space).  If $\ba\in\AZD$, and $\dU\subset\ZD$,
then we define $\ba_\dU\in\sA^\dU$ by
$\ba_\dU:=[a_\fu]_{\fu\in\dU}$.  If $\fz\in\ZD$, then strictly speaking,
$\ba_{\fz+\dU}\in\sA^{\fz+\dU}$; however, it is sometimes convenient
to `abuse notation' and treat $\ba_{\fz+\dU}$ as an element of
$\sA^{\dU}$ in the obvious way.  If $\dX\subset\dY\subseteq\ZD$, and
$\bx\in\sA^\dX$ and $\by\in\sA^\dY$, we 
write ``$\bx \sqsubset \by$'' if $\bx=\by_{\dX}$.

\parag{Cellular automata:} For any $\fv\in\ZD$,
we define the {\dfn shift} $\shift{\fv}:\AZD\into\AZD$ by
$\shift{\fv}(\ba)_{\fz} = a_{\fz+\fv}$ for all $\ba\in\AZD$ and
$\fz\in\ZD$.  A {\dfn cellular automaton} is a
transformation $\Phi:\AZD\into\AZD$ that is
continuous and commutes with all shifts.
Equivalently, $\Phi$ is determined by
a {\dfn local rule} $\phi:\sA^\Nh\into\sA$ so that
$\Phi(\ba)_{\fz} = \phi(\ba_{\fz+\Nh})$ for all $\ba\in\AZD$ and
$\fz\in\ZD$ \cite{Hedlund}.   Here, $\Nh\subset\ZD$ is a finite
set which we normally imagine as a `neighbourhood of the origin'.
If $\Nh\subseteq \dB(r):=\CC{-r...r}^D$, we 
say that $\Phi$ has {\dfn radius} $r$.

\parag{Subshifts:}
  A subset $\gA\subset\AZD$ is a {\dfn subshift}
\cite{LindMarcus,Kitchens} if $\gA$ is closed in the Cantor topology,
and if $\shift{\fz}(\gA)=\gA$ for all $\fz\in\ZD$.  For any
$\dU\subset\ZD$, we define $\gA_\dU:=\set{\ba_\dU}{\ba\in\gA}$. 
In particular, for any $r>0$, let 
$\Adm := \gA_{\dB(r)}$ be the set of {\dfn admissible
$r$-blocks} for $\gA$.  We say $\gA$ is {\dfn subshift of
finite type} (SFT) if there is some $r>0$ (the {\dfn radius} of $\gA$)
such that  $\gA$ is entirely described by $\Adm$, in the sense
that $\gA=\set{\ba\in\AZD}{\ba_{\dB(\fz,r)}\in\Adm, \
\forall\fz\in\ZD}$ (here, $\dB(\fz,r):=\fz+\CC{-r...r}^D$).  If $D=1$, then 
a {\dfn Markov subshift} is an SFT $\gA\subset\sA^\Zahl$ determined
by a set $\gA_{\{0,1\}}\subset\sA^{\{0,1\}}$ of {\dfn admissible transitions}; 
equivalently, $\gA$ is the set of all bi-infinite directed paths in a 
digraph  whose vertices are the elements of $\gA$, with
an edge $a\leadsto b$ iff $(a,b)\in\gA_{\{0,1\}}$. 
If $D=2$, then let $\dE_1:=\{(0,0),(1,0)\}$ and
 $\dE_2:=\{(0,0),(0,1)\}$.  A {\dfn Wang subshift}
 is an SFT $\gA\subset\AZD[2]$
determined by sets $\gA_{\dE_1}\subset\sA^{\dE_1}$ and
$\gA_{\dE_2}\subset\sA^{\dE_2}$ of {\dfn edge-matching conditions}.
Equivalently, $\gA$ is the set of all {\dfn tilings} of the plane $\RD[2]$
by unit square tiles (corresponding to the elements of $\sA$) with
notched edges representing the edge-matching conditions 
\cite[Ch.11]{GrunbaumShephard}.
 
If $\sX$ is any set and $F:\gA\into\sX$ is a function, then $F$ is
{\dfn locally determined} if there is some {\dfn radius} $r\in\Natur$
and some {\dfn local rule} $f:\Adm\into\sX$ such that
$F(\ba)=f(\ba_{\dB(r)})$ for any $\ba\in\gA$.  If $\sX$ is any
discrete space, then $F:\gA\into\sX$ is continuous iff $F$ is locally
determined.  For example, if $\sA$ and $\sB$ are finite sets, then a
(subshift) {\dfn homomorphism} is a continuous, $\shift{}$-commuting
function $\Phi:\BZD\into\AZD$ (e.g. a CA is a homomorphism with
$\sA=\sB$); it follows that $F(\ba)=\Phi(\ba)_0$ is locally
determined.  If $\gB\subset\BZD$ is a subshift of finite type, and
$\Psi:\BZD\into\AZD$ is a homomorphism, then $\gA:=\Psi(\gB)\subset\AZD$
is called a {\dfn sofic shift}.  If $D=1$, we define the {\dfn
language} of a subshift $\gA\subset\sA^\Zahl$ by
$\sL(\gA):=\Union_{n=0}^\oo \gA_\CC{0...n}$; then $\gA$ is sofic if and
only if $\sL(\gA)$ is a {\dfn regular language}; i.e. a language
recognized by a finite automaton \cite[\S2.5 \& \S9.1]{HopcroftUllman}.
Equivalently, $\gA$ is the set of all bi-infinite directed paths in a
digraph whose vertices are (nonbijectively) $\sA$-labelled.

\ignore{
  If $\gA\subset\AZD$ is a subshift and $r>0$, then a word
$\ba\in\sA^\dU$ is {\dfn $r$-admissible} if $\ba_{\fu+\dB(r)} \in
\Adm$ whenever $\fu+\dB(r)\subset\dU$.  If $\gA$ is a SFT of radius
$r$, then $\ba$ is called {\dfn admissible} if $\ba$ is
$r$-admissible.  For any $R>r>0$, we say that $\gA$ is {\dfn
$(R,r)$-extensible} if, for any $\bb\in\sA^{\dB(R)}$, if $\bb$ is
$r$-admissible then $\ba\in\Adm[R]$ ---i.e.  there is some $\ba\in\gA$ so
that $\bb=\ba_{\dB(R)}$.  (For example, any subshift is trivially
$(R,R)$-admissible for every $R>0$).  If $\gA$ is an SFT of radius
$r$, then $\gA$ is {\dfn extensible} (or has the {\dfn extension
property} \cite[Dfn.2.1]{Sch98}) if there is some $R_0>0$ so that
$\gA$ is $(R,r)$-extensible for all $R\geq R_0$.  (Conversely, if
$\gA$ is any shift which is $(R,r)$-extensible for all sufficiently
large $R$, then $\gA$ must be an SFT of radius $r$, by a compactness
argument.)
}

  If $\Phi:\AZD\into\AZD$ is a cellular automaton, then we say $\gA$
is  {\dfn $\Phi$-invariant} if $\Phi(\gA)=\gA$, and $\gA$ is
{\dfn weakly $\Phi$-invariant} if $\Phi(\gA)\subseteq\gA$
(i.e. $\Phi$ is an {\dfn endomorphism} of $\gA$).  For example:

 (a)\ The set $\Fix{\Phi}:=\set{\ba\in\AZD}{\Phi(\ba)=\ba}$  of
$\Phi$-{\dfn fixed points} is a $\Phi$-invariant
SFT.  Likewise, if $p\in\Natur$ and $\fv\in\ZD$, then the set
$\Fix{\Phi^p}$ of  {\dfn $(\Phi,p)$-periodic points} and the set
$\Fix{\Phi^p\circ\shift{-p\fv}}$ of {\dfn $(\Phi,p,\fv)$-travelling waves} are
$\Phi$-invariant SFTs.

  (b)\ If $\dP\subset\ZD$ is a finite-index subgroup, then the set  
$\sA^{\ZD/\dP}:=\Intsct_{\fp\in\dP} \Fix{\shift{\fp}}$
of {\dfn $\dP$-periodic configurations}
is a weakly $\Phi$-invariant SFT.
(Furthermore, all $\Phi$-orbits in $\sA^{\ZD/\dP}$ are eventually
periodic, because $\sA^{\ZD/\dP}$ is finite.)  

  (c) If $\Phi$ has radius $R$, then for any $r>0$,
$\Phi$ induces a function $\Phi:\sA^{\dB(R+r)}\into\sA^{\dB(r)}$.
If $\gA\subset\AZD$ is an SFT determined by a set $\gA_{(r)}\subset\sA^{\dB(r)}$ of
admissible $r$-blocks, then $\statement{$\Phi(\gA)\subseteq\gA)$}
\iff \statement{$\Phi(\gA_{(r+R)}) \subseteq \gA_{(r)}$}$. 
Furthermore, if $D=1$ and $\gA\subset\AZ$ is irreducible, then 
$\Phi(\gA)=\gA$ (surjectively) iff
$\Phi$ is finite-to-one on $\gA$
iff $\Phi$ has no `diamonds' i.e. words $\ba,\bc\in\gA_{(r)}$
and $\bb\neq\bb'\in\gA_L$ ($L\geq R$)
with $\Phi(\ba\bb\bc)=\Phi(\ba\bb'\bc)$ 
\cite[Thm.8.1.16]{LindMarcus}.
The endomorphism set of an SFT is quite huge;
see \cite[Ch.3]{Kitchens} or
\cite[\S13.2]{LindMarcus}.

\ignore{  For example, if
$\gA\subset\AZ$ is mixing, then $\Aut{\gA,\shift{}}$ contains copies
of $\Zahl^{\dirsum \oo}$, the free group on two generators, and any
finite group.  }

 (d)\ Let $\Phi^\oo(\AZD):=\Intsct_{t=1}^\oo \Phi^t(\AZD)$ be the
{\dfn eventual image} of $\Phi$.  Then $\Phi^\oo(\AZD)$ is a
$\Phi$-invariant subshift (possibly non-sofic), which contains
$\Fix{\Phi^p\circ\shift{-p\fv}}$ for any $p\in\Natur$ and $\fv\in\ZD$.

\ignore{
  (e)  $\gA\subset\AZ$ is called {\dfn rank one} if, for any $\eps>0$,
there is some `tile' $\bt\in\sA^H$ (for some $H\in\Natur$)
such that any $\ba\in\gA$ can be covered with copies of $\bt$
except for a set of \Cesaro density $\eps$.  The endomorphisms of
$\gA$ can then be `approximated' by shifts \cite{KingWkClos}.  Special
cases of this result apply to Toeplitz shifts
(\cite{BK1,BK2,Downarowicz} and \cite[\S4]{DKL}) and quasisturmian
shifts \cite[\S14]{PivatoQuasi}.}

\parag{Trails:}
If $\fy,\fz\in\ZD$, then we write
``$\fy\adjacent\fz$'' if $|\fz-\fy|= 1$.   A {\dfn trail} is a sequence
$\zeta=(\fz_1\adjacent\fz_2\adjacent\cdots\adjacent\fz_n)$.  A subset
$\dY\subset\ZD$ is {\dfn trail-connected} if, for any $\fx,\fy\in\dY$,
there is a trail $\fx=\fz_0\adjacent\fz_1\adjacent
\cdots\adjacent\fz_n=\fy$ in $\dY$.

\parag{Font conventions:}
Upper case calligraphic letters ($\sA,\sB,\sC,\ldots$) denote 
finite alphabets.  Upper-case Gothic letters ($\gA,\gB,\gC,\ldots$)
are subsets of $\AZD$ (e.g. subshifts), lowercase bold-faced letters
($\ba,\bb,\bc,\ldots$) are elements of $\AZD$, and Roman letters
($a,b,c,\ldots$) are elements of $\sA$ or ordinary numbers.
Lower-case sans-serif ($\ldots,\fx,\fy,\fz$) are elements of $\ZD$,
upper-case hollow font ($\dU,\dV,\dW,\ldots$) are subsets of $\ZD$,
upper-case Greek letters ($\Phi,\Psi,\ldots$) denote functions
on $\AZD$ (e.g. CA), and lower-case Greek letters ($\phi,\psi,\ldots$)
denote other functions (e.g. local rules).

\section{Domains and Defects\label{S:defect}}

  Let $\gA\subset\AZD$ be any subshift.  If $\ba\in\AZD$, then the
{\dfn defect field} $\energy:\ZD\into\Natur\union\{\oo\}$ is defined
\[
\forall\ \fz\in\ZD, \quad
  \energy(\fz)\quad:=\quad \max\set{r\in\Natur}{\ba_{\dB(\fz,r)} \in \gA_{(r)}}.
\]
It is easy to see that $\energy$ is `Lipschitz' in the sense that
$|\energy(\fy)-\energy(\fz)|\leq|\fy-\fz|$.  
The {\dfn defect set} of $\ba$ is the set $\dD(\ba)\subset\ZD$ of local minima of
$\energy$.  

\EXAMPLE{\label{X:defect.energy}
(a) $\ba\in\gA$ if and only if $\energy(\fz)=\oo$ for some (and thus all) $\fz\in\ZD$.  In this case, $\dD(\ba)=\emptyset$. 

(b) Suppose $\gA$ is an SFT determined by a set $\gA_{(r)}\subset\sA^{\dB(r)}$
of admissible $r$-blocks, and
let $\dX:=\set{\fz\in\ZD}{\ba_{\dB(\fz,r)}\not\in\gA_{(r)}}$.
Assume for simplicity
that $\gA_{(r-1)}=\sA^{\dB(r-1)}$.  Then 
$\energy(\fz)=r+d(\fz,\dX)$, where
$d(\fz,\dX):=\D \min_{\fx\in\dX} \ |\fz-\fx|$.
In particular, $\energy(\fz)=r$ if and only if $\fz\in\dX$,
and this is the smallest possible value for $\energy(\fz)$.
Thus, $\dD(\ba)=\dX$.

(c) Let $\sA:=\{0,1\}$.
Let $\gS$ be the sofic shift defined by the $\sA$-labelled digraph 
$\textcircled{\footnotesize 1} \leftrightarrows
\textcircled{\footnotesize 0}  \leftrightarrows
\textcircled{\footnotesize 0}$  (this is
the invariant sofic shift of ECA\#18 
mentioned in the introduction).  Let
\[
\bs\quad:=\quad
 [\ldots 00\ 01 \ 00\ 01 \ 01 \ \underbrace{00 \ 00 \ 00 \ 00  \overbrace{00}^\dX  00 \
00 \ 00 \ 00}_\dY \ 10 \ 00 \ 10 \ 00 \ 00 \ 10 \ldots].
\]
Then for any $\fz\in\Zahl$, $\energy[\bs](\fz)$ is the distance
from $\fz$ to the {\em farthest} endpoint of $\dY$;  this is the maximum
radius $r$ such that $\bs_\CC{z-r...z+r}$ is $\gS$-inadmissible.  Thus,
$\energy[\bs](\fz)$ takes a minimal value of 8 (inside $\dX$) and
increases linearly in either direction.  Thus, $\dD(\bs)=\dX$.

(d)  Let $\sA=\sB\union\sD$ (think of $\sB$ and $\sD$ as `invariant
subalphabets' in the sense of \cite{Elo93a,Elo93b}),
and let $\gA:=\sB^{\ZD}\union\sD^{\ZD}$.  Let $\sC:=\sB\intsct\sD$ be
the (possibly empty) set of `ambiguous symbols', 
and let $\sB^*:=\sB\setminus\sC$ and $\sD^*:=\sD\setminus\sC$.
Any $\ba\in\AZD$
is a mixture of $\sB^*$-symbols, $\sC$-symbols, and $\sD^*$-symbols.
If $\fz\in\ZD$ and $a_\fz\in\sB^*$, then
$\energy(\fz)=\min\set{|\fy-\fz|}{a_\fy\in\sD^*}$.  
If $a_\fz\in\sD^*$, then $\energy(\fz)=\min\set{|\fy-\fz|}{a_\fy\in\sB^*}$.
If $a_\fz\in\sC$, then $\energy(\fz)=
\min\{r \ ; \ a_\fx\in\sB^* \AND a_\fy\in\sB^* $ $  \ \mbox{for some 
$\fx,\fy\in\dB(\fz,r)$}\}$.  Thus, $\dD(\ba)$ is the set of all points 
which are either on a `boundary' between a $\sB^*$-domain and a $\sD^*$-domain,
or which lie roughly in the middle of a $\sC$-domain.
}

\Remarks{If $\gA$ is an SFT [Example \ref{X:defect.energy}(b)], then
the set $\dD(\ba)$ encodes all information about the defects of
$\ba$. However, if $\gA$ is {\em not} of finite type, then $\dD(\ba)$
is an inadequate description of the larger-scale `defect structures'
of $\ba$.  Thus, instead of treating the defect as a precisely defined
subset of $\ZD$, it is better to think of it as a `fuzzy' object
residing in the low areas in the defect field $\energy$.  This is particularly
appropriate when we want to evaluate a locally determined function
(such as a cellular automaton or eigenfunction) on $\ba$ at a point
`close' (but not {\em too} close) to the defect.  The advantage of
this definition is its applicability to any kind of subshift (finite
type, sofic, or otherwise).  Nevertheless, most of our examples will
be SFTs, and we may then refer to the specific region $\dD\subset\ZD$
as `the defect'. This is the approach taken in \cite{PivatoDefect0},
for example. }

Let
\[
 \tlgA\quad:=\quad\set{\ba\in\AZD}{\D\sup_{\fz\in\ZD} \ \energy(\fz) \ = \ \oo}
\]
be the set of `slightly defective' configurations.   
If $\ba\in\tlgA\setminus\gA$, then we say $\ba$ is {\dfn defective}.
Elements of $\tlgA$ may have infinitely large defects, but also have
arbitrarily large non-defective regions.
Clearly $\gA\subset\tlgA$, and $\tlgA$ is a $\shift{}$-invariant,
dense subset of $\AZD$ (but not a subshift).

For any $R>0$, let
$\unflawed[R](\ba):=\set{\fz\in\ZD}{\energy(\fz)\geq R}$.  Thus,
$\ba\in\tlgA$ iff $\unflawed[R](\ba)\neq\emptyset$ for all $R>0$.  For
example, if $\gA$ is an SFT determined by a set $\gA_{(r)}$ of admissible
$r$-blocks, and $\dD=\set{\fz\in\ZD}{\ba_{\dB(\fz,r)}\not\in\gA_{(r)}}$ as
in Example \ref{X:defect.energy}(b), then
\[
\unflawed[R](\ba)\quad=\quad\set{\fz\in\ZD}{d(\fz,\dD)\geq R-r}
\quad=\quad\ZD\setminus\dB(\dD,R-r).
\]
  However, if $\gA$ is {\em not} an SFT,
then in general $\unflawed[R](\ba) \neq\ZD\setminus\dB(\dD,R')$ for
any $R'>0$.

\Proposition{\label{defect.dimension}}
{
 Let $\Phi:\AZD\into\AZD$ be a CA with radius $r>0$.
\bthmlist
  \item Let $\gA\subset\AZD$ be a weakly $\Phi$-invariant subshift. 
Then $\Phi(\tlgA)\subseteq\tlgA$.

  \item If $\ba\in\tlgA$, and $\ba'=\Phi(\ba)$,
then $\energy[\ba'] \geq \energy[\ba]-r$.
Thus, for all $R\in\Natur$, \ 
 $\unflawed[R+r](\ba)\subseteq \unflawed[R](\ba')$.
\ethmlist
}
\bthmprf  {\bf(b)} Let $\fz\in\ZD$ and suppose $\energy(\fz)=R$.  Thus,
$\ba_{\dB(\fz,R)}\in\gA_R$.  But $\gA$ is $\Phi$-invariant;
hence $\ba'_{\dB(\fz,R-r)}\in\gA_{(R-r)}$.  Hence $\energy[\ba'](\fz)\geq R-r$.
Then {\bf(a)} follows from {\bf(b)}.
\ethmprf

If $\ba\in\tlgA$, we say $\ba$ has a {\dfn range $r$ domain boundary}
 if $\unflawed(\ba)$ is trail-disconnected. Domain boundaries divide
 $\ZD$ into different `domains', which may correspond to different
 transitive components of $\gA$ (see \S\ref{S:interface}), different
 eigenfunction phases (see \S\ref{S:disl}), or different cocycle
 asymptotics \cite[\S2.3]{PivatoDefect2}.

\begin{figure}[h]
\centerline{\footnotesize
\begin{tabular}{cc}
\includegraphics[scale=0.3]{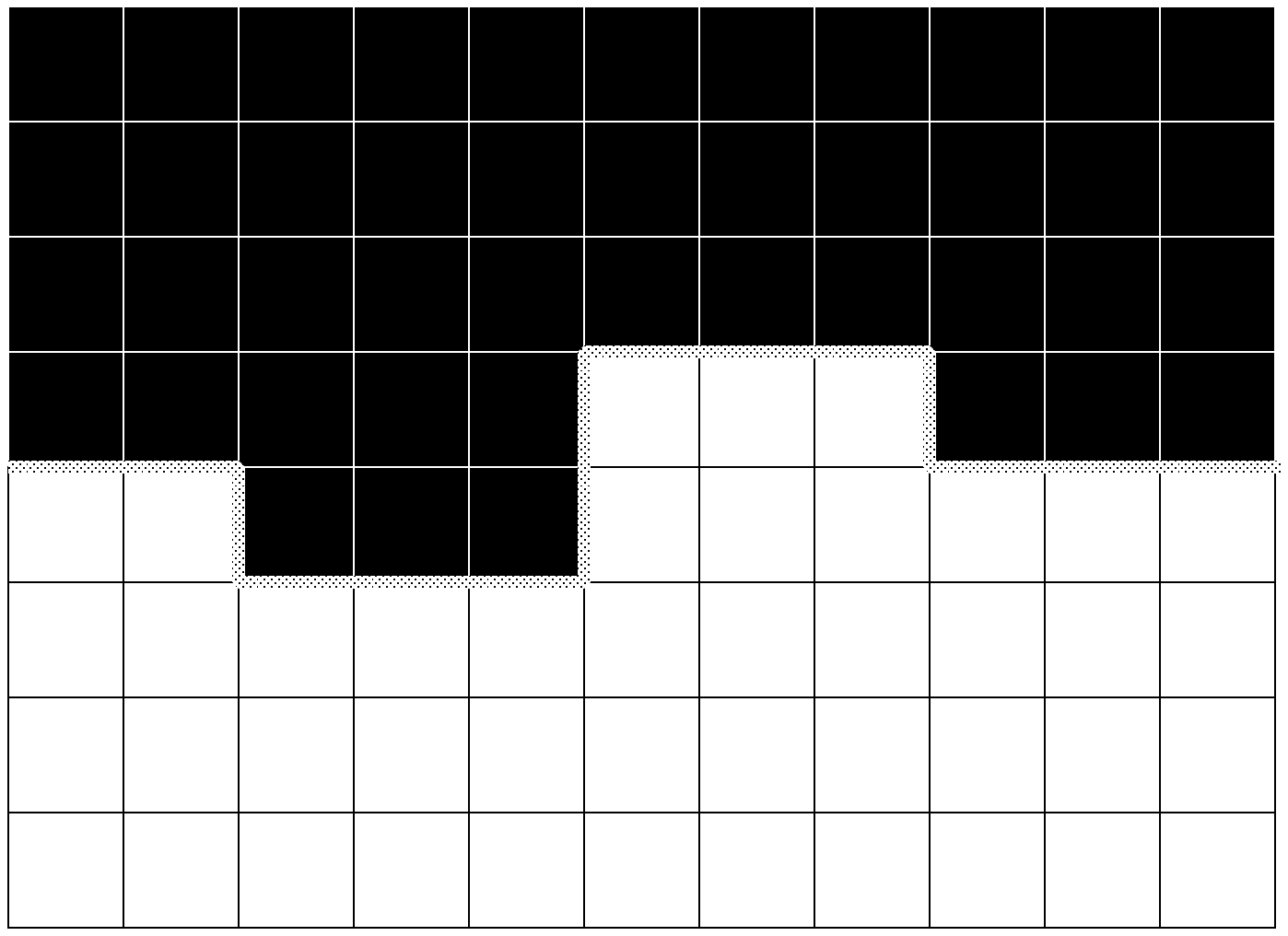} &
\includegraphics[scale=0.3]{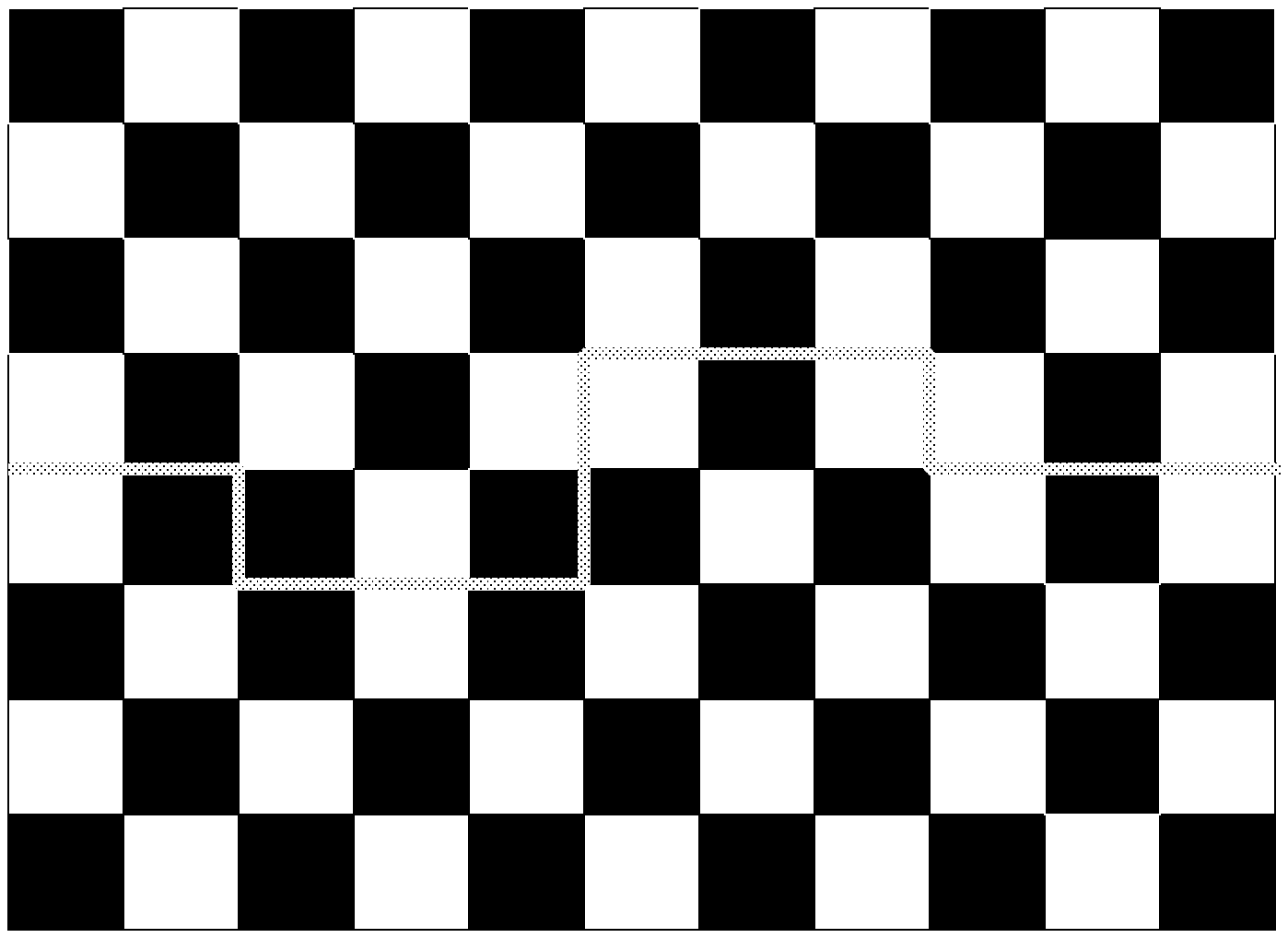} \\
{\bf(A)} & {\bf(B)} \\
\includegraphics[scale=0.5]{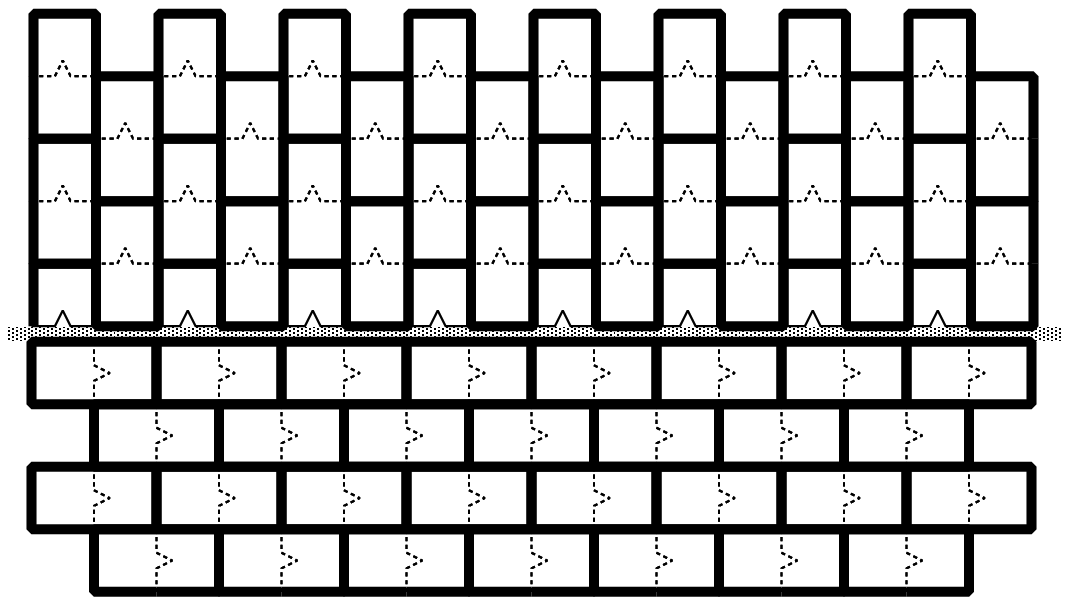}&
\includegraphics[scale=0.5]{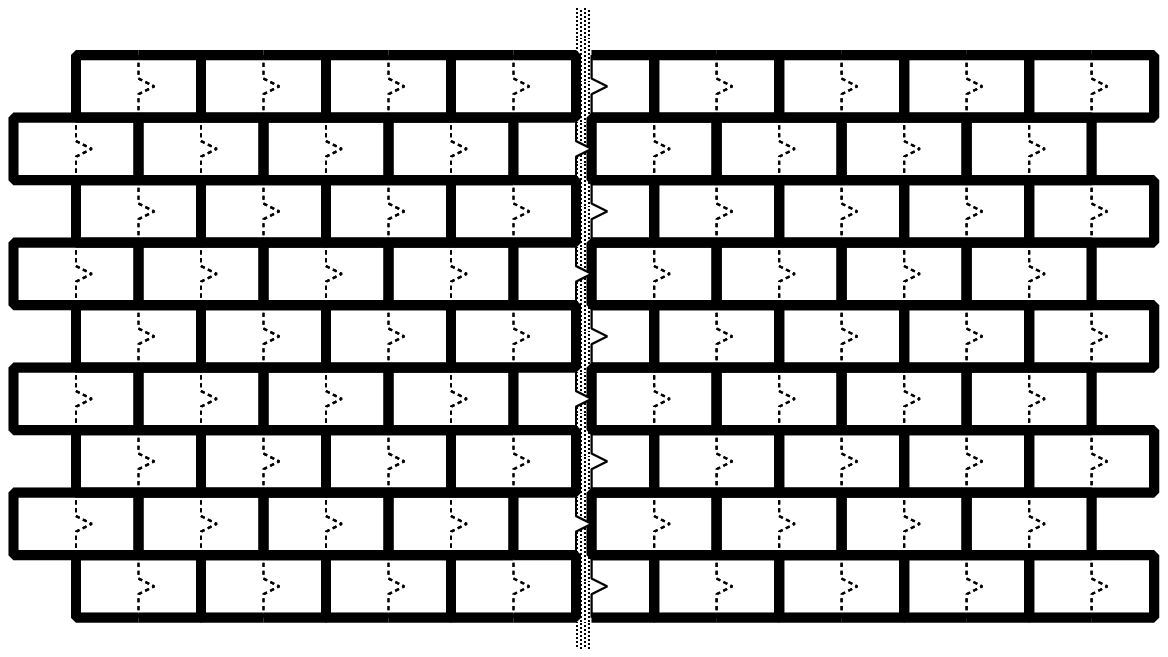}\\
{\bf(C)}  & {\bf(D)}  
\end{tabular}}
\Caption{fig:codim}{Domain boundaries.
{\bf(A)} \ An interface in $\Mono$;
see Examples \ref{X:codim.one}(a) and \ref{X:interface}(b).
{\bf(B)} \ A dislocation in $\Checker$;
see Examples \ref{X:codim.one}(b) and  \ref{X:dislocation}(e).
{\bf(C,D)} \  Boundaries in $\Dom$;
 see Examples \ref{X:codim.one}(c) and
 \ref{X:interface}(c)\projdisl{ and \ref{X:irrational.dislocation}(c)}.}
 \end{figure}

\EXAMPLE{\label{X:codim.one}(a) 
Let $\sA=\{\black,\white\}$, and let
$\Mono\subset\AZD[2]$ be the {\dfn monochromatic} SFT defined by the condition that
no $\black$ can be adjacent to a $\white$.
Figure \ref{fig:codim}(A) shows a domain boundary in $\Mono$.  See also Example \ref{X:interface}(b).

 (b)  
Let $\sA=\{\black,\white\}$, and let
$\Checker\subset\AZD[2]$ be the {\dfn checkerboard} SFT defined by the condition that
no $\black$ can be adjacent to a $\black$,
and no $\white$  can be adjacent to a $\white$.
Figure \ref{fig:codim}(B) shows a domain boundary in $\Checker$.   See also Example \ref{X:dislocation}(e).

  (c)  Let $\sD := \lb\{\raisebox{-0.5em}{
\includegraphics[scale=0.65]{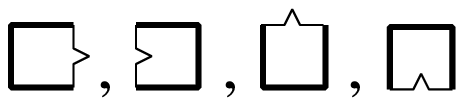}}\rb\}$,
 and let $\Dom\subset\sD^{\ZD[2]}$ 
be the {\dfn domino} SFT defined by the obvious edge-matching conditions.
Figure \ref{fig:codim}(D,E)  shows two domain boundaries in $\Dom$;
see also Example\projdisl{s} \ref{X:interface}(c)\projdisl{ and \ref{X:irrational.dislocation}(c)}.

  (d)  Let $\sA=\{0,1\}$, and let  $\gS\subset\AZD$ and $\bs\in\tl\gS$
be as in Example \ref{X:defect.energy}(c).
Then $\bs$ has a domain boundary in the region $\dY$.
Figure \ref{fig:defect.intro}(A) shows defects of this
type evolving under the iteration of ECA\#18.
\projdisl{See also Example \ref{X:irrational.dislocation}(a).}
}

  Domain boundaries are defects of {\em codimension one} (they
disconnect the space), and are the only kind which exist in
one-dimensional cellular automata.  In two-dimensional cellular
automata, there are also defects of codimension {\em two}, which do
{\em not} disconnect the space, but instead resemble localized
`holes'.  In three-dimensional CA, a codimension-two defect has the
topology of an extended `string', while a `hole'-shaped defect has
codimension {\em three}.  The precise definition of defect codimension
involves homotopy groups; see \cite{PivatoDefect2}.  The present paper
is concerned only with defects of codimension one ---i.e. domain
boundaries.

\ignore{A connected component $\dY$
of $\unflawed(\ba)$ is called {\dfn projective} if $\diam{\dY}=\oo$
(where $\diam{\dY}:=\sup\set{|\fy_1-\fy_2|}{\fy_1,\fy_2\in\dY}$). }

\parag{Projective domain boundaries:}

The action of a cellular automaton may locally change the geometry of
a defect, and we are mainly interested in properties that are
invariant under such local modifications. 
  If $\ba\in\tlgA$, then a connected component $\dY$ of 
$\unflawed(\ba)$ is called {\dfn projective} if, for all
$R\geq r$, \ $\dY\intsct\unflawed[R](\ba)\neq \emptyset$.
We say that $\ba\in\tlgA$ has a {\dfn projective domain boundary} if there is
some $R\geq0$ such that $\unflawed[R](\ba)$ has at least two
projective components.  (Hence $\unflawed(\ba)$ is disconnected for
all $r\geq R$.)

 A trail-connected subset $\dY\subset\ZD$ is
{\dfn spacious} if, for any $R>0$, there exists $\fy\in\dY$ with
$\dB(\fy,R)\subset\dY$. 
If $\dY$ is a connected component of $\unflawed(\ba)$, 
then ($\dY$ is projective) $\IMPLIES$ ($\dY$ is spacious).  If $\gA$ is
a subshift of finite type, then the converse is also true.

\EXAMPLE{\label{X:projective.spacious} (a) A proper subset
$\dY\subset\Zahl$ is spacious iff either $\dY=\OC{-\oo...Z}$ or
$\dY=\CO{Z...\oo}$ for some $Z\in\Zahl$.

(b) We say the defect in $\ba$ is {\dfn finite} if
 $\D \lim_{|\fz|\goto\oo}\ \energy(\fz)=\oo$.
(This implies that $\dD(\ba)$ is finite. If $\gA$ is a subshift of
finite type, then the converse is also true.)
If $\gA\subset\AZ$, then $\ba\in\tlgA$ has a projective domain boundary iff 
the defect of $\ba$ is finite.
If $R$ is large enough, then 
 $\unflawed[R](\ba)=\OC{-\oo...Z_1}\disj \CO{Z_2...\oo}$
for some $Z_1,Z_2\in\Zahl$.
}

\parag{Essential vs. Persistent Defects:}
 Let $\Phi:\AZD\into\AZD$ be a cellular automaton, and suppose that
$\phi(\gA)=\gA$.  If $\ba\in\tlgA$ then $\ba$ has a {\dfn
$\Phi$-persistent} defect if, for all $t\in\Natur$, \
$\ba'=\Phi^t(\ba)$ is also defective.  Otherwise $\ba$ has a {\dfn
transient} defect ---i.e. one which eventually disappears.  Our main 
goal is to determine when defects are persistent.  We say $\ba$ has a
{\dfn removable} defect if there is some $r>0$ and some $\ba'\in\gA$
so that $a'_\fz=a_\fz$ for all $\fz\in\unflawed(\ba)$ (i.e. the defect
can be erased by modifying $\ba$ in a finite radius of the defective
region).  Otherwise $\ba$ has an {\dfn essential} defect.

\EXAMPLE{\label{X:removable.finite.defects} Let $\gA\subset\AZ$
be a subshift of finite type.
Then  $(\gA,\shift{})$ is topologically mixing if and only if no finite
 defect is essential.  }

\Proposition{\label{biject.essential.persistent}}
{
Let $\Phi:\AZD\into\AZD$ be a CA and let
$\gA\subset\AZD$ be a $\Phi$-invariant subshift.
If \linebreak $\Phi:\gA\into\gA$ is bijective, then any essential defect is 
$\Phi$-persistent.
}
\bthmprf
Suppose $\ba\in\tlgA$ has an essential defect and let $\ba':=\Phi(\ba)$.
We must show that $\ba'$ is also defective.  We will suppose not and derive
a contradiction.
Suppose that $\Phi$ has radius $H>0$.
Then for any $r>0$, we have a map $\Phi_r:\gA_{(r+H)}\into\gA_{(r)}$.

\Claim{There exists $R>0$ such that, for all $r\geq R$, the function
$\Phi_r:\gA_{(r+H)}\into\gA_{(r)}$ is bijective.}
\bclaimprf
  Suppose not.  Let $\{r_n\}_{n=1}^\oo$ be a sequence such that,
for each $n\in\Natur$, $\Phi_{r_n}:\gA_{(r_n+H)}\into\gA_{(r_n)}$ is not bijective.
Let $\bc_n\in\gA_{(r_n)}$ be a point with two $\Phi_{r_n}$-preimages
in $\gA_{(r_n+H)}$, say $\bb_n$ and $\bb'_n$.  By dropping to a subsequence
if necessary, we can assume $\bc_1\sqsubset\bc_2\sqsubset\cdots$
and $\bb_1\sqsubset\bb_2\sqsubset\cdots$ and $\bb'_1\sqsubset\bb'_2\sqsubset\cdots$.  Thus, there are limit points
$\bc,\bb,\bb'\in\gA$ such that $\bc_n=\bc_{\dB(r_n)}$, 
$\bb_n=\bb_{\dB(r_n)}$, and $\bb'_n=\bb'_{\dB(r_n)}$ for all $n\in\Natur$
(because $\gA$ is compact).  Also, $\Phi(\bb)=\bc=\Phi(\bb')$
(because $\Phi$ is continuous). But $\bb\neq\bb'$ (by construction) so
this contradicts the supposed bijectivity of $\Phi$.
\eclaimprf
If $\ba'$ is not defective, then $\ba'\in\gA$.  Let $\bb\in\gA$ be
the unique $\Phi$-preimage of $\ba'$ in $\gA$ (recall that
$\Phi:\gA\into\gA$ is bijective).  Let $R>0$ be
as in Claim 1, and let $\dY:=\unflawed[R+H](\ba)$.  

\Claim{$\ba_{\dY}=\bb_{\dY}$.}
\bclaimprf
  If $\fy\in\dY$, then $\ba_{\dB(\fy,R+H)}\in\gA_{(R+H)}$,
so $\Phi(\ba_{\dB(\fy,R+H)})
\ = \ (\ba')_{\dB(\fy,R)}$ is an element of $\gA_{(R)}$.  But 
by definition, 
$\Phi(\bb_{\dB(\fy,R+H)}) \ = \ \ba'_{\dB(\fy,R)}$ also.
Thus, $\ba_{\dB(\fy,R+H)} = \bb_{\dB(\fy,R+H)}$
(because  Claim 1 says $\Phi_r:\gA_{(r+H)}\into\gA_{(r)}$ is bijective.)
This is true for every $\fy\in\dY$.  
\eclaimprf
Thus, $\ba$ and $\bb$ are equal on $\unflawed[R+H](\ba)$,
so the defect in $\ba$ is removable.  Contradiction.
\ethmprf
\Corollary{}
{
 Suppose $\gA\subseteq\Fix{\Phi}$ or 
 $\gA\subseteq\Fix{\Phi^p\circ\shift{p\fv}}$ 
{\rm(for some $p\in\Natur$ and $\fv\in\ZD$)}.  If $\ba\in\tlgA$, 
then any essential defect in $\ba$ is $\Phi$-persistent.\qed}

\begin{figure}
\centerline{
\begin{tabular}{cc}
\includegraphics[clip=true, trim=14 50 0 0, scale=2]{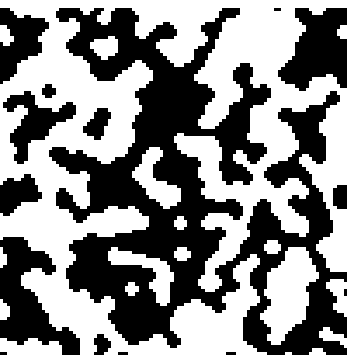}
&
\includegraphics[clip=true, trim=14 50 0 0, scale=2]{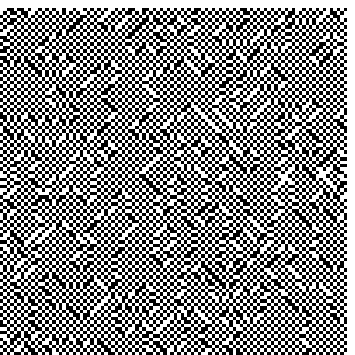}\\
{\bf(A)} & {\bf(B)} \\
\end{tabular}
}
\Caption{fig:voter}{Configurations with many (nonprojective)
domain boundaries.  {\bf(A)} \ An equilibrium of a  voter
CA [Example \ref{X:biject.persistent.defect}(a)], with $\Mono$-domain boundaries.
\ 
{\bf(B)} \ An equilibrium of the zero-temperature antiferromagnet 
 [Example \ref{X:biject.persistent.defect}(b)], with $\Checker$-domain boundaries.}
\end{figure}

\EXAMPLE{\label{X:biject.persistent.defect}
(a) Let $\Mono$ be as in Example \ref{X:codim.one}(a).
The domain boundary in Figure \ref{fig:codim}(A)
is essential because it separates two infinite
domains of opposite colour (and one infinite domain would have to be
completely erased to eliminate the defect).  Let $\Phi:\AZD[2]\into\AZD[2]$ be a
{\dfn Voter CA} \cite{Vichniac} with local rule $\phi:\sA^\Nh\into\sA$ defined
\[ 
  \phi(\ba) \quad:=\quad \choice{\black &\If& N(\ba) < \theta;\\
\white &\If& N(\ba) \geq \theta;}
\quad
\mbox{where}\quad
 N(\ba)\  := \  \frac{\#\set{\nh\in\Nh}{a_\nh=\white}}{\#(\Nh)},
\]
and where $\theta\in\CC{0,1}$ is some {\dfn threshold}.
Then $\Mono\subset\Fix{\Phi}$;  hence the domain boundary in Figure \ref{fig:codim}(A) is $\Phi$-persistent.  If $\theta$ is close to 1/2
(for example, in the {\dfn zero-temperature ferromagnet} CA, $\theta:=1/2$), then
domain boundaries like Figure \ref{fig:codim}(A) are roughly stationary.
The CA rapidly evolves from random initial conditions
to a mottled equilibrium configuration with infinitely many such
boundaries, as in Figure \ref{fig:voter}(A) (these boundaries are not
essential, because the domains are finite).
If $\theta$ is close to $0$ (respectively to $1$), then the boundary in Figure \ref{fig:codim}(A) will
rapidly propagate north (respectively south).  Likewise,
a small `seed' of one colour will
grow monotonically like a crystal. (The asymptotic shape of these `crystals'
has been studied in \cite{Wi78,Gr93,Gr98b,Gr98a}.)

(b) Let $\Checker$ be as in
Example \ref{X:codim.one}(b).
The domain boundary in Figure \ref{fig:codim}(B)
is essential because it separates
two infinite domains of opposite `phase'.
 Let $\Phi\colon\AZD[2]\goto\AZD[2]$ be a
{\dfn zero-temperature antiferromagnet} CA,
with local rule $\phi:\sA^\Nh\into\sA$ defined
\beq 
  \phi(\ba) &:=& \choice{\black &\If& N_0(\ba)-N_1(\ba) <0;\\
\white &\If& N_0(\ba)-N_1(\ba) \geq 0;}\\
\mbox{where}
\quad 
N_0(\ba) &:=& \#\set{\nh=(h_0,h_1)\in\Nh}{a_\nh=\white \AND h_0+h_1=0 \bmod{2}}, \\
\AND N_1(\ba) &:=& \#\set{\nh=(h_0,h_1)\in\Nh}{a_\nh=\white \AND h_0+h_1=1 \bmod{2}}. 
\eeq
Then $\Checker\subset\Fix{\Phi}$; 
hence the domain boundary in Figure \ref{fig:codim}(B) is $\Phi$-persistent.
Figure \ref{fig:voter}(B) shows a $\Phi$-fixed configuration with
many $\Checker$-domain boundaries.
}

\section{\label{S:interface}Interfaces}

 A $\Phi$-invariant subshift  $\Sft\subset\AZD$ is
{\dfn $(\Phi,\shift{})$-transitive} if, for any nonempty
open subset $\gO\subset\gA$, the union 
\[\Union_{t\in\Natur} \ \Union_{\fz\in\ZD} 
\Phi^{-t}\shift{-\fz}(\gO)\]
 is dense in $\gA$.  Equivalently, there
exists  $\ba\in\gA$ such that the orbit
$\{\Phi^t\circ\shift{\fz}(\ba)\}_{t\in\Natur, \fz\in\ZD}$ is dense
in $\gA$.

 Let $\Sft\subset\AZD$ be a $(\Phi,\shift{})$-nontransitive subshift,
and suppose $\gA=\gA_1\disj\cdots\disj\gA_K$, where
$\gA_1,\ldots,\gA_K$ are clopen $(\Phi,\shift{})$-transitive
components.  As $\gA_1,\ldots,\gA_K$ are clopen, their indicator
functions are locally determined; hence there is some $r>0$, and a
function $\kap:\gA_{(r)}\into\CC{1...K}$ such that, for any
$\ba\in\gA$, \ we have $\lb(\ba\in\gA_k\rb)
\Leftrightarrow\lb(\kap(\ba_{\dB(r)})=k\rb)$.  For any $\fz\in\ZD$,
let $\kap_\fz(\ba) \ := \ \kap(\ba_{\dB(\fz,r)})$.  Hence, if
$\ba\in\Sft_k$, then $\kap_\fz(\ba) = k$ for all $\fz\in\ZD$. However,
$\kap_\fz(\ba)$ is also well-defined for any $\ba\in\tlgA$ such that
$\ba_{\dB(\fz,r)}$ is $\Sft$-admissible.  Hence, if $\ba\in\tlgA$,
then $\kap_\fz(\ba)$ is well-defined for all $\fz\in\unflawed(\ba)$.
If $\fy,\fz\in\unflawed(\ba)$, then we say that $\ba$ has an {\dfn $(\gA,\Phi)$-interface}
between $\fy$ and $\fz$ if $\kap_\fy(\ba) \neq \kap_{\fz}(\ba)$.
(We will simply call this an {\dfn interface} when $\gA$ and $\Phi$ are clear
from context).

\Proposition{\label{interface.prop}}
{
Let $\Phi:\AZD\into\AZD$ be a CA and let $\gA\subset\AZD$ be a subshift
with $\Phi(\Sft)\subseteq\Sft$.
Suppose $\ba\in\tlgA$ has an $(\gA,\Phi)$-interface.
Let $\kap$ have radius $r'$ and let $r\geq r'+1$.  Then
\bthmlist 
  \item The interface in $\ba$ is a range $r$ domain boundary.

  \item $\kap(\ba)$ is constant on each connected component of 
$\unflawed(\ba)$.
 Suppose $\unflawed(\ba)$ has connected components
 $\{\dY_n\}_{n=1}^N$ {\rm(where $N\in\Natur\union\{\oo\}$)}.
There is a function $\sK:\CC{1...N}\into\CC{1...K}$ such that,
for any \linebreak $n\in\CC{1...N}$ and any $\fy\in\dY_n$,
$\kap_\fy(\ba)=\sK(n)$.

\item If two of the components $\dY_n$ and $\dY_m$ are projective,
and $\sK(n)\neq\sK(m)$, then there is an essential domain boundary
between $\dY_n$ and $\dY_m$.
\ethmlist
}
\bthmprf {\bf(b)}  Suppose $\fx$ and $\fz$ are in the same connected
component of $\unflawed(\ba)$.  Let $\fx=\fy_0\adjacent\fy_1\adjacent\ldots\adjacent\fy_N=\fz$ be a trail from $\fx$ to $\fz$ in $\unflawed(\ba)$.

\Claim{For all $n\in\CC{1...N}$, \ $\kap_{\fy_n}(\ba)=\kap_{\fy_{n-1}}(\ba)$.}
\bclaimprf
If $\dU:=\dB(\fy_{n},r)$ then $\ba_\dU\in\gA_\dU$
(because $\energy(\fy_n)\geq r$ by definition). 
Thus, there exists $\bb\in\gA$ such that $\bb_\dU=\ba_\dU$.  
There is a (unique) $k\in\CC{1...K}$ such that $\bb\in\gA_k$. 
Thus, $\kap_{\fy_{n-1}}(\ba) \ \eeequals{(*)} \ \kap_{\fy_{n-1}}(\bb)=k=\kap_{\fy_n}(\bb)\ \eeequals{(\dagger)} \ \kap_{\fy_n}(\ba)$,
where $(*)$ is because $\dB(\fy_{n-1},r')\subset\dU$ and 
$(\dagger)$ is because $\dB(\fy_{n},r')\subset\dU$.
\eclaimprf
Apply Claim 1 inductively to conclude that $\kap_\fx(\ba)=\kap_{\fz}(\ba)$.

{\bf(a)} If there
exist $\fy,\fz\in\unflawed(\ba)$ with
$\kap_\fy(\ba) \neq \kap_{\fz}(\ba)$, then  {\bf(b)} says
$\fy$ and $\fz$ must be in different connected components of $\unflawed(\ba)$.
\quad {\bf(c)}  follows immediately.
\ethmprf
\ignore{The {\dfn rank} of the interface is the number $N$ of connected components
of $\unflawed(\ba)$.
The function $\sK:\CC{1...N}\into\CC{1...K}$ is the {\dfn signature}
of the interface.  ($\sK$ is only defined up the
(arbitrary) ordering of $\dY_1,\ldots,\dY_N$.)
The {\dfn projective rank}  is the number $M\leq N$ of
projective components of $\unflawed(\ba)$
(hence $M>1$ iff $\ba$ has a {\em projective} domain boundary).}

By reordering if necessary, assume that the projective components are
$\dY_1,\ldots,\dY_M$. 
We call the restricted function $\sK_*:\CC{1...M}\into\CC{1...K}$
the {\dfn signature} of the interface.
Proposition \ref{interface.prop}(c) says the interface is essential if
$\sK_*$ is not constant.

\begin{figure}[h]
\centerline{\footnotesize
\begin{tabular}{|c|c|c|c|c|c|c|c|}
\hline
\includegraphics[scale=3]{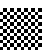} &
\includegraphics[scale=3]{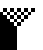} &
\includegraphics[scale=3]{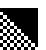} &
\includegraphics[scale=3]{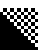} &
\includegraphics[scale=3]{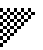} &
\includegraphics[scale=3]{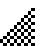} &
\includegraphics[scale=3]{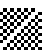} &
\includegraphics[scale=3]{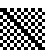}\\
$(*)$  &
$\boldsymbol{(\beta)}$  & 
$\boldsymbol{(\alpha^+)}$ & $\boldsymbol{(\omg^+)}$ &
 $\boldsymbol{(\alp^-)}$ & $\boldsymbol{(\omg^-)}$
& $\boldsymbol{(\gamma^-)}$ & $\boldsymbol{(\gamma^+)}$ \\
\hline
\end{tabular} }
\Caption{fig:184.defects}{$(*)$  The periodic background generated
by $\ECA{184}$ acting on  $\gG$.\quad
$\boldsymbol{(\alpha^\pm,\beta,\omg^\pm)}$:  Interfaces in
$\gG$; see Example \ref{X:persistent.interface}(a).
$\boldsymbol{(\gamma^\pm)}$: Dislocations  in $\gG$;
see Examples \ref{X:persistent.dislocation}(c)
and \ref{X:permanent.dislocations}(c).
(In the nomenclature of \cite[\S III(A)]{BNR91},
$\gam^+=\vec{0}_2$, \  $\gam^-=\stackrel{\leftarrow}{1}_2$, 
while both $\alp^+$ and $\omg^+$ are ``$\vec{0}_\oo$'', and both $\alp^-$ and $\omg^-$ are  ``$\stackrel{\leftarrow}{1}_\oo$'').
See also \cite[Example 1.2(a)]{PivatoDefect0}}.
\end{figure}

\EXAMPLE{\label{X:interface}
(a) \label{X:persistent.interface}
(ECA\#184)  Let $\sA=\{\black,\white\}$, and let $\gG=\gG_0\disj\gG_1\disj\gG_*$,
where $\gG_0:=\{\overline{\black}\}$, $\gG_1:=\{\overline{\white}\}$, 
and $\gG_*:=\{\overline{\black\white}, \ \overline{\white\black}\}$. 
(Here, $\overline{\black}:=[\ldots\black\black.\black\black\ldots]$ and 
$\overline{\black\white}:=[\ldots\black\white.\black\white\black\white\ldots]$, etc., where the period is before the zeroth coordinate). 
If $\ECA{184}$ is ECA\#184, then $\gG_0\union\gG_1 \subset
\Fix{\ECA{184}}$, while
$\ECA{184}\restr{\gG_*} = \shift{}$,
as shown in Figure \ref{fig:184.defects}$(*)$.
Thus $\gG$ has three $(\ECA{184},\shift{})$-transitive components, so there are six
possible interfaces:
{\footnotesize \[\begin{array}{lcr}
\begin{array}{rc}
\alp^+: & \gG_* \
 \fbox{$\ldots\black\white\black\white\black\white$}
\ \fbox{$\black\black\black\black\black\black\ldots$} \ \gG_0 \\
\omg^+: & \gG_0 \
  \fbox{$\ldots\black\black\black\black\black\black$} \
\fbox{$\black\white\black\white\black\white\ldots$} \  \gG_* \\
\bet: &  \gG_0 \ \fbox{$\ldots\black\black\black\black\black\black$} \
\fbox{$\white\white\white\white\white\white\ldots$} \ \gG_1 
\end{array}
&&
\begin{array}{rc}
\alp^-: & \gG_* \
  \fbox{$\ldots\black\white\black\white\black\white$} 
\ \fbox{$\white\white\white\white\white\white\ldots$}  \ \gG_1 \\
\omg^-: & \gG_1 \
 \fbox{$\ldots\white\white\white\white\white\white$} 
\  \fbox{$\black\white\black\white\black\white\ldots$} \ \gG_* \\
\eps: &  \gG_0 \ \fbox{$\ldots\white\white\white\white\white\white$}
\ \fbox{$\black\black\black\black\black\black\ldots$} \ \gG_1 
\end{array}
\end{array}\]}
\ignore{(Here, we use the convention that $\black=0$ and $\white=1$, to ease comparison between equations and figures.)}
Figure \ref{fig:184.defects}$(\alp^\pm,\bet,\omg^\pm)$
shows the $\ECA{184}$-evolution of these defects. (The $\eps$ defect
is unstable, and immedately `decays' into  $\omg^-$ and
$\alp^+$ defects travelling in opposite directions.)
Figure \ref{fig:defect.intro}(D) showed the long-term
$\ECA{184}$-evolution of these defects.
If $\bg\in\tl\gG$ has a finite defect, then $\bg$ can be
written as an ensemble of range-$r$ defects $\bd_1,\ldots,\bd_N$
 arranged along a line, 
with $\dY_0,\ldots,\dY_N$ being the $\gG$-admissible intervals between
these defects:

\centerline{$
 \cdots-\!\!-\!\dY_0\!\longrightarrow \ \bd_1
\longleftarrow\!\dY_1\!\longrightarrow \ \bd_2 \
\longleftarrow\!\dY_2\!\longrightarrow \cdots
\longleftarrow\!\dY_{N-1}\!\longrightarrow \ \bd_N \
\longleftarrow\!\dY_N\!-\!\!- \cdots$}

The  projective  components are $\dY_0$ and $\dY_N$.
Hence the interface is essential if $\sK(0)\neq\sK(N)$.

(b)  Let $\sA=\{\black,\white\}$ and let $\Mono\subset\AZD[2]$ be
as in  Example \ref{X:codim.one}(a).
Then $\Mono=\{\black^\oo,\white^\oo\}$, where $\black^\oo$ is the
solid black configuration, and 
$\white^\oo$ is the solid black configuration.
Let $\Phi:\AZD[2]\into\AZD[2]$ be an CA such that $\Mono\subset\Fix{\Phi}$
[e.g. a voter CA from Example \ref{X:biject.persistent.defect}(a)]
 Then
$\Mono$ has two $(\Phi,\shift{})$-transitive components, $\gM_0:=\{\black^\oo\}$ and  
$\gM_1:=\{\white^\oo\}$.
Figure \ref{fig:codim}(A) shows a domain boundary in $\Mono$.
If $\dY_0$ is the northern connected component and $\dY_1$ is
the southern component, then we have $\sK(0)=0$ and $\sK(1)=1$.
Both components are projective, so this is an essential interface.

(c)   Let $\Dom$ be as in Example \ref{X:codim.one}(c).
Despite appearances, 
the domain boundary in Figure \ref{fig:codim}(C) is {\em not} an
interface, because $(\Dom,\shift{})$ is topologically transitive
\cite[Lemma 2.1]{Ein}.  Instead, this is 
a `gap' defect; see \cite[Example 2.14(b)]{PivatoDefect2}.
}
The next result implies that
the defects in Examples \ref{X:interface}(a,b) must be persistent.

\Proposition{\label{persistent.interface}}
{
Let $\Phi:\AZD\into\AZD$ be a CA. If $\gA\subset\AZD$ is a subshift
with $\Phi(\Sft)=\Sft$, then any essential $(\gA,\Phi)$-interface is
$\Phi$-persistent.
If $\ba\in\tlgA$ has an essential interface, then
$\Phi(\ba)$ also has an essential interface, with the same 
signature as $\ba$.
}
\bthmprf Let $\ba':=\Phi(\ba)$ and suppose  $\Phi$ has radius $R>0$.
Then each projective component of $\unflawed[r+R](\ba)$ is
contained in a projective component of $\unflawed(\ba')$,
because $\unflawed[r+R](\ba)\subseteq\unflawed[r](\ba')$
by Lemma \ref{defect.dimension}(b).
Let $\unflawed[R+r](\ba)$  have projective components
$\dY_1,\ldots,\dY_M$, and let
$\unflawed(\ba')$  have projective components
$\dY'_1,\ldots,\dY'_M$, where $\dY_m\subset\dY'_m$ for all $m\in\CC{1...M}$.
If $\fy\in\dY_m\subset\dY'_m$, then 
$\kap_\fy(\ba)$ and $\kap_\fy(\ba')$ are well-defined,
and, in the notation of  Proposition \ref{interface.prop}(b),
we must have $\kap_\fy(\ba')=\sK(m)=\kap_\fy(\ba)$. 
\ethmprf

\noindent{\em Remark \refstepcounter{thm}\label{transitive.decomposition.remark}\thethm:}  
  If $\gA$ is not $\shift{}$-transitive, then a {\dfn
$\shift{}$-transitive decomposition} of $\gA$ is a collection of
disjoint clopen $\shift{}$-transitive components $\gA_1,\ldots,\gA_R$
such that $\gA=\gA_1\disj\cdots\disj\gA_R$.  (Not all
non-transitive subshifts admit such a decomposition.)  If $\gA$ has a
$\shift{}$-transitive decomposition, and $\Phi(\gA)=\gA$, then $\Phi$
induces a permutation $\varphi:\{\gA_1,...,\gA_R\} \into
\{\gA_1,...,\gA_R\}$, and each $(\Phi,\shift{})$-transitive component
of $\gA$ is a union of all elements of $\{\gA_1,...,\gA_R\}$ in some
$\varsigma$-orbit.  (In particular, if $\varphi=\Id{}$, then the
$(\Phi,\shift{})$-transitive components of $\gA$ are also
$\gA_1,\ldots,\gA_R$.)

%\subsection{Rational Dislocations \label{S:dislocation}}
\section{\label{S:disl}Dislocations}

  In a periodic crystalline solid, a {\em dislocation}
(or {\em fault line}) is an internal surface separating two domains whose
crystal structures are spatially out of phase.  We will use the word
{\dfn dislocation} to describe an analogous domain boundary in
a configuration in $\AZD$.  
The main results of this section
are Theorems \ref{gap.prop}, \ref{persistent.dislocation}
 \ref{ae.spec}, and \ref{persistent.proj.disl}.

\EXAMPLE{\label{X:one.dim.per.SFT}
If $D=1$ and $\gA\subset\AZ$ is a $\shift{}$-transitive SFT, then any
dislocation in $\gA$ takes a simple form.  For simplicity,
suppose $\gA$ is a  Markov subshift.
For any $a,c\in\sA$, say that $c$ is {\dfn reachable from $a$
in time $t$} if there is a word $\bb\in\sA^{t-1}$ such that $a\bb c$ is
$\gA$-admissible.  There is some
$P=P(\gA)\in\Natur$ (the {\dfn period} of $\gA$) and a {\dfn phase partition}
$\sA=\sA_0\disj\cdots\disj\sA_{P-1}$ such that, if $a\in\sA_n$ and
$c\in\sA_m$, then $c$ is reachable from $a$ in time $t$ only if $t
\equiv m-n$ (mod $P$) \cite[Prop.4.5.6]{LindMarcus}.  Hence,
$\gA$ is mixing iff $P=1$.  A sequence $\ba=[\ldots a_{-1} a_0 a_1
\ldots]\in\AZ$ thus has a {\dfn dislocation} at $0$ if $a_0\in\sA_n$
and $a_1 \in\sA_{m}$, but $m\neq n+1 \pmod{P}$.  The {\dfn phase
gap} of the dislocation is the value $m-(n+1)$ (as an element of
$\Zahlmod{P}$).  Two such defects can `cancel out' if and only if
their phase gaps together sum to zero, mod $P$.

  For example, let $\sA:=\{a,b,c\}$ and let $\gA\subset\AZ$
be the Markov subshift defined by the $\sA$-labelled digraph
$\textcircled{\footnotesize $a$} \leftrightarrows
\textcircled{\footnotesize $b$}  \leftrightarrows
\textcircled{\footnotesize $c$}$.  Then $P=2$, with
$\gA_0=\{a,c\}$ and $\gA_1=\{b\}$.  Hence,
the sequence $[\ldots ababcbabc.abcbcbab\ldots]$ has a dislocation at the 
decimal point.
}

  If $\gA\subset\AZ$ is a non-finite type subshift, or if
$\gA\subset\AZD$ is of finite type, for $D\geq 2$, then dislocations
and their `phase gaps' can take a more subtle form than in Example 
\ref{X:one.dim.per.SFT}.  Spectral
theory provides the tools to characterize these
dislocations. 

\subsection{\label{S:dislocation}Rational Dislocations}

 Let $\sC(\Sft)$ be the $\Cplx$-vector space of continuous $\Cplx$-valued
functions on $\Sft$.
Let $\dT\subset\Cplx$ be the unit circle.  A {\dfn $(\Phi,\shift{})$-eigenfunction}
of $\Sft$ is any $f\in\sC(\Sft)$ admitting some {\dfn generalized eigenvalue}
 $\blam = (\lam_0; \ \lam_1,\ldots,\lam_D)\in\Torus{D+1}$
 such that:  
\bdesc
\item[{\bf(a)}]  $f\circ \Phi  =  \lam_0 f$.
\item[{\bf(b)}] For any $\fz=(z_1,\ldots,z_D)\in\ZD$,
\ $f\circ \shift{\fz}  = \blam^{\fz} f$, 
where $\blam^{\fz} := \lam_1^{z_1}\cdots \lam_D^{z_D}$. 
\edesc
 Let $\Spec\subset\Torus{D+1}$ be the set of all
such eigenvalues.
 For any $\blam\in\Spec$, let 
\[
\Eigenspace \quad:= \quad\Eigenspace(\gA,\Phi,\shift{}) \quad=\quad 
\set{f\in\sC(\Sft)}{f\circ \Phi  =  \lam_0 f \AND
f\circ\shift{\fz}=\blam^\fz, \ \forall \fz\in\ZD}
\]
be the {\dfn eigenspace} of $\blam$.  
We next relate $\Spec[\gA,\shift{}]$ to $\Spec$,
and review basic spectral theory.

\Lemma{\label{spectral.lemma}}
{
  Let $\Phi\colon\AZD\rightarrow\AZD$ be a CA. 
Let $\gA\subset\AZD$ be a subshift with $\Phi(\gA)=\gA$.
\bthmlist
  \item $\Spec$ is a multiplicative subgroup of $\Torus{D+1}$.
Let $\blam_1,\blam_2\in\Spec$ with
$f_1\in\Eigenspace[\blam_1]$ and $f_2\in\Eigenspace[\blam_2]$.
If $\blam=\blam_1\cdot\blam_2$ and $f=f_1\cdot f_2$, then 
$f\in\Eigenspace[\blam]$.

  \item  If $\SFT$ is transitive, then
 $\dim(\Eigenspace)=1$ for each $\blam\in\Spec$.

  \item Suppose $\Psi:\BZD\into\BZD$ is another CA,  and
$\gB\subset\BZD$ is a $\Psi$-invariant subshift.  Let 
\linebreak$\xi:\SFT\into(\gB,\Psi,\shift{})$ be an epimorphism.  Then 
$\Spec[\gB]\subseteq\Spec[\gA]$.
For any $\blam\in \Spec[\gB]$,
there is a linear monomorphism $\xi_*:\Eigenspace(\gB)\into\Eigenspace(\gA)$ defined by $\xi_*(f) = f\circ \xi$.

\item If $\gA$ is $\shift{}$-transitive,  then there is a homomorphism
$\tau:\Spec[\gA,\shift{}]\into\dT$ such that
\[
\Spec \quad  = \quad   \set{(\lam_0;\blam)}{\blam\in\Spec[\gA,\shift{}]
 \AND  \lam_0=\tau(\blam)} \quad \cong \quad  \Spec[\gA,\shift{}].
\]
{\rm[For example, if $\Phi\restr{\gA}\equiv\shift{\fz}$, then
$\tau(\blam)=\blam^\fz$,  for all $\blam\in\Spec[\gA,\shift{}]$.]}

Furthermore, if $\blam\in \Spec[\gA,\shift{}]$ and
$\lam_0:=\tau(\blam)$, then
 $\Eigenspace[(\lam_0;\blam)](\gA,\Phi,\shift{})
  =  \Eigenspace(\gA,\shift{})$.

\item  Suppose $\gA$ is not  $\shift{}$-transitive, but has
$\shift{}$-transitive decomposition $\D \gA=\Disj_{n=1}^N\gA_n$.
Then:
\bitem
\item[{\rm[i]}] $\Phi$ induces a permutation $\varphi:\{\gA_1,...,\gA_K\}
\into \{\gA_1,...,\gA_K\}$, and $(\gA,\Phi,\shift{})$ is transitive
iff $\varphi$ is transitive.  In this case,
$\Spec[\gA_1,\shift{}] = \cdots = \Spec[\gA_K,\shift{}]$.

\item[{\rm[ii]}] There is then a homomorphism
$\tau:\Spec[\gA_1,\shift{}]\into\dT$ such that
\ignore{
$\Spec \ = \ 
\set{(\rho\cdot \tau(\blam); \ \blam)}{
\blam\in\Spec[\gA_1,\shift{}] \ 
\& \ \mbox{$\rho\in\dT$ a $K$th root of unity}}$

\hspace{5.5em} $\cong \ \Zahlmod{K}\x\Spec[\gA,\shift{}]$.
}
\beq
\Spec & = & 
\set{(\rho\cdot \tau(\blam); \ \blam)}{
\blam\in\Spec[\gA_1,\shift{}] \ \mbox{and $\rho\in\dT$ is a $K$th root of unity}}\\
&\cong & \Zahlmod{K}\x\Spec[\gA,\shift{}].
\eeq
\eitem

\item Let $\blam\in\Spec$.  The following are equivalent:
\bitem
  \item[{\rm[i]}] The subgroup $\set{\blam^\fz}{\fz\in\ZD}\subset\dT$ is finite.

  \item[{\rm[ii]}] $\blam = (\lam_0,\ldots,\lam_D)$, where
$\lam_0,\ldots,\lam_D$ are complex roots of unity.

  \item[{\rm[iii]}]  Every  $F\in\Eigenspace$ 
is locally determined (hence $F(\gA)\subset\Cplx$ is finite).\qed
\eitem
\ethmlist
}
\bthmprf
{\bf(a,b,c,f):} See \cite[\S1.5]{Fogg}, \cite[\S5.5]{Walters},
or \cite[Prop.2.53]{KurkaBook}.
%(\cite[Thm.3.1]{Walters} or \cite[Thm.4.2]{Petersen} for measurable version):

\begin{sloppypar}
{\bf(d)} For each $\blam\in\Spec[\gA,\shift{}]$,  
part {\bf(c)} yields a linear map $\Phi_*\colon\Eigenspace(\gA,\shift{})
 \into  \Eigenspace(\gA,\shift{})$. Part {\bf(b)} says
that $\dim[\Eigenspace(\gA,\shift{})]=1$, so
there exists $\lam_0\in\Cplx$
such that $f\circ\Phi = \lam_0 f$ for all $f\in \Eigenspace(\gA,\shift{})$.
But  $|f\circ\Phi|=|f|$, hence $|\lam_0|=1$, so $\lam_0\in\dT$.
Define $\tau(\blam):=\lam_0$.  It follows that $f\in \Eigenspace[(\lam_0;\blam)]$.  Use part {\bf(a)} to check that $\tau$ is a homomorphism.
\end{sloppypar}

{\bf(e)} [i] follows from Remark \ref{transitive.decomposition.remark}
and part {\bf(c)}.

[ii]  By reordering if necessary,
assume $\varphi(\gA_k)=\gA_{k-1}$ for all $k\in\CC{1...K}$ (where
$k-1$ is mod $K$).
For any $\blam\in\Spec[\gA_k,\shift{}]$, part {\bf(c)} yields a monomorphism
$\Phi_*\colon\Eigenspace(\gA_{k},\shift{})\rightarrow
\Eigenspace(\gA_{k+1},\shift{})$.  Let $\blam\in\Spec[\gA_1,\shift{}]$,
let  $f_1\in\Eigenspace(\gA_1,\shift{})$, and let 
$f_k := f_1\circ\Phi^{k-1}$, for all  $k\in\CC{2...K}$.
Then  $f_k \in \Eigenspace(\gA_{k},\shift{})$, and $\{f_1,\ldots,f_K\}$
is a $\Cplx$-basis for $\D \Eigenspace(\gA,\shift{}) = \Dirsum_{k=1}^K \Eigenspace(\gA_{k},\shift{})$.   Furthermore,
$\Phi_*^K:\Eigenspace(\gA_{1},\shift{})\into\Eigenspace(\gA_{1},\shift{})$ is
a linear map of a one-dimensional space; hence, just as in {\bf(d)}, 
there is some $\tau=\tau(\blam)\in\dT$ such that $f_1\circ\Phi^K = \tau^K\cdot f_1$.  Thus, $f_k\circ\Phi^K = \tau^K\cdot f_k$ for all $k\in\CC{1...K}$, and thus, $F\circ\Phi^K=\tau^K F$ for every $F\in \Eigenspace(\gA,\shift{})$
(because $\Eigenspace(\gA,\shift{})$ is spanned by $\{f_1,\ldots,f_K\}$).

``$\supseteq$''
Let $\rho$ be a $K$th root of unity and let $\lam_0:= \rho\cdot \tau$.
If $F := \D \sum_{k=1}^K \lam_0^{K-k} f_k$, then $F\in \Eigenspace(\gA,\shift{})$,
and  $F\circ\Phi = \lam_0 F$, so  $F\in \Eigenspace[\lam_0](\gA,\Phi)$;
hence $F\in\Eigenspace[(\lam_0;\blam)](\gA,\Phi,\shift{})$ as desired.

``$\subseteq$'' If $\lam_0\in\dT$ and $F\in\Eigenspace[(\lam_0;\blam)](\gA,\Phi,\shift{})$,
then $\lam_0^KF \ = \ F\circ\Phi^K \ = \ \tau^K F$ (where $\tau=\tau(\blam)$),
so $\lam_0^K = \tau^K$, so
 $\lam_0=\rho\tau$, for some $K$th root of unity $\rho$.
\ignore{Also, $F'\in \Eigenspace(\gA,\shift{})$, so $F'=\sum_{k=1}^K c_k f_k$,
where $\{f_1,\ldots,f_k\}$ are as above, and $c_1,\ldots,c_K\in\Cplx$.
But then $c_k = \lam c_{k+1}$ for all $k\in\CC{1...K}$ (where $c_{K+1}=c_1$),
because 
\[
\sum_{k=1}^K \lam c_k f_k
\lam F' \ = \
 F'\circ\Phi \ = \ \sum_{k=1}^K c_k f_k\circ\Phi \ = \ 
\sum_{k=1}^K c_k f_{k+1} \ = \ \sum_{k=1}^K c_{k-1} f_{k}.
\]
Hence, we conclude that $c_k = \lam_0^{K-k} c_1$, which means
$F'=c_1 F$, where $F$ is as above. }
\ethmprf

We say that $\blam$ is {\dfn rational} if any (and thus all) of the
conditions in Lemma \ref{spectral.lemma}{(f)} hold.  If $\SFT$ is
transitive, then we define the {\dfn radius} of $\blam$ to be the
radius of any nontrivial $F\in\Eigenspace$ [this is finite by Lemma
\ref{spectral.lemma}{(f)}[iii], and independent of $F$ by Lemma
\ref{spectral.lemma}{(b)}].  Let $\RatSpec\subset\Spec$ be the
subgroup of rational eigenvalues.
Lemma \ref{spectral.lemma}{(d,e)} implies that
$\RatSpec$ is nontrivial iff $\RatSpec[\gA,\shift{}]$ is nontrivial.
Meanwhile, $\RatSpec[\gA,\shift{}]$ is nontrivial iff
$(\gA,\shift{})$ has a periodic factor.

\ignore{, and
if $\gA$ is $\shift{}$-transitive, then 
Lemma \ref{spectral.lemma}{(d)} allows us to 
blur the distinction between $\RatSpec[\gA,\shift{}]$  and
$\RatSpec[\gA,\Phi,\shift{}]$ (and the corresponding eigenspaces).}

Let
$\hRatSpec:=\{\mbox{continuous homomorphisms} \
\del\colon\RatSpec\rightarrow\dT\}$ be the {dual group}
of $\RatSpec$.  For any $(t;\fz)\in\Zahl\x\ZD$, define
$\displace{(t;\fz)}\in\hRatSpec$ by $\displace{(t;\fz)}(\blam) := \lam_0^t
\lam_1^{z_1}\cdots \lam_D^{z_D}$.  We will see that $\displace{(t;\fz)}$
corresponds to a `displacement' in time by $t$ and in space by $\fz$.
 The group homomorphism $\bdel:\Zahl\x\ZD\ni (t;\fz)\mapsto
\displace{(t;\fz)}\in\hRatSpec$ has dense image
$\hRatSpec$ (and in most of our examples, is surjective).
Hence we will regard elements of $\hRatSpec$ as `generalized
space-time shifts', and call them {\dfn displacements}.
If $\gA$ is $\shift{}$-transitive, then 
Lemma \ref{spectral.lemma}{(d)} yields a natural isomorphism
$\hRatSpec\cong\hRatSpec[\gA,\shift{}]$,
so all displacements are `space shifts'. 

\EXAMPLE{\label{X:ratspec}
(a) Let $\gA\subset\AZ$ be a  one-dimensional $\shift{}$-transitive SFT, with
period $P$ and phase partition $\sA_0\disj\cdots\disj\sA_{P-1}$,
as in Example \ref{X:one.dim.per.SFT}. If $\lam:=e^{2\pi\bi/P}$, then 
\[
\Spec[\gA,\shift{}]\quad=\quad\RatSpec[\gA,\shift{}]
\quad=\quad \{\lam^q\}_{q=0}^{P-1}\quad\cong\quad\Zahlmod{P}.
\]
(the group of $P$th roots of unity).
For example,  define $f:\gA\into\dT$ by 
$f(\ba):=\lam^q$ iff $a_0\in\sA_q$.
Then $f\in\Eigenspace[\lam](\gA,\shift{})$.  

  If $\Phi(\gA)=\gA$, then there is {\dfn phase rotation}
$r\in\Zahlmod{P}$ such that, for any $\ba\in\gA$, if $a_0\in\sA_p$,
then $\Phi(\ba)_0 \in\sA_{p+r}$.  The homomorphism
$\tau:\Spec[\gA,\shift{}]\into\dT$ in Lemma
\ref{spectral.lemma}{(d)} is then defined $\tau(\lam^q):=\lam^{rq}$ for
all $q\in\Zahlmod{P}$.  To see this, note that $f\circ\Phi = \lam^{rq}
f$ for any $f\in\Eigenspace[\lam^q](\gA,\shift{})$; hence
$\Eigenspace[\lam^q](\gA,\shift{}) =
\Eigenspace[(\lam^{rq};\lam^q)](\gA,\Phi,\shift{})$.  In this case,
$\hRatSpec \cong \hRatSpec[\gA,\shift{}]\cong \Zahlmod{P}$, and the group
homomorphism  $\bdel:\Zahl\ni \fz\mapsto\displace{(0,\fz)}\in \hRatSpec$ is
surjective, with kernel $P\Zahl$.

(b) Let $\gA\subset\AZD$ be $\shift{}$-transitive,
and let $\dP\subset\ZD$ be a finite-index
subgroup. We say that $\Sft$ is {\dfn $\dP$-periodic} if
$\gA\subset\Fix{\shift{\fp}}$ for all $\fp\in\dP$.
 Let $\tlfZ:=\ZD/\dP$ be the quotient group, with quotient map
$\ZD\ni \fz\mapsto \tlfz\in\tlfZ$. 
Let $\fe_1,\ldots,\fe_D\in\ZD$ be the unit vectors.  Then
$\{\tlfe_1,\ldots,\tlfe_D\}$ generates $\tlfZ$.
For each $d\in\CC{1...D}$, let $p_d\in\Natur$ be the (finite) order
of $\tlfe_d$ in $\tlfZ$, and let 
$\lam_d:=e^{2\pi\bi/p_d}$.  Then
\[
\Spec[\gA,\shift{}]\quad=\quad
\RatSpec[\gA,\shift{}]\quad=\quad \set{(\lam_1^{n_1},...,
\lam_D^{n_D})}{n_1,\ldots,n_D\in\Zahl}.
\]  
 [For example, fix $\ba_0\in\gA$, and $\forall \ \tlfz\in\tlfZ$, let
 $\ba_{\tlfz}:=\shift{\fz}(\ba_0)$ (well-defined because
 $\ba_0$ is $\dP$-periodic).  Then
 $\gA=\{\ba_{\tlfz}\}_{\tlfz\in\tlfZ}$ because $\gA$ is
 $\shift{}$-transitive. Let $\blam:=(\lam_1^{n_1},...,\lam_D^{n_D})$,
and define $f\colon\Sft\rightarrow\dT$ by $f(\ba_{\tlfz}) 
 :=  \blam^\fz, \ \forall \ \fz\in\ZD$.  Then
 $f\in\Eigenspace(\gA,\shift{})$.]
The  homomorphism $\bdel:\ZD\ni \fz\mapsto\displace{(0,\fz)}\in \hRatSpec[\gA,\shift{}]$ is surjective with
kernel $\dP$, so $\hRatSpec[\gA,\shift{}]\cong\tlfZ$. 
 If $\Phi(\gA)=\gA$, then there exists $\tlfz\in\tlfZ$ such that
$\Phi\restr{\gA}=\shift{\fz}$. The function 
$\tau:\Spec[\gA,\shift{}]\into\dT$ in Lemma \ref{spectral.lemma}(d)
is then defined by $\tau(\blam):=\blam^\fz$.}  

\Remarks{The homomorphism $\bdel:\Zahl\into\hRatSpec[\gA,\shift{}]$ is not always surjective.  For example, if $p\in\Natur$, and
 $\gA$ is a $p$-adic Toeplitz shift \cite{Downarowicz} or a 
nonperiodic substitution shift
induced by a substitution $\sA\into\sA^p$ \cite[Thm.7.3.1]{Fogg}, then
$\hRatSpec\cong\Zahl(p)$ is the $p$-adic integers, and $\bdel$ is the 
natural embedding $\Zahl \hookrightarrow\Zahl(p)$, which is {\em not} surjective.}

  Suppose $\SFT$ is transitive.  Let $\blam\in\RatSpec$ have radius $r>0$,
 and let $f\in\Eigenspace$.  Then for
any $\ba\in\Sft$ we can write $f(\ba) = f(\ba_{\dB(r)})$.
  For any $\fz\in\ZD$, let $f_\fz(\ba) \ := \
f(\ba_{\dB(\fz,r)})$.  Hence, if $\ba\in\Sft$, then $f_\fz(\ba) =
f(\shift{\fz}(\ba)) \ = \ \blam^\fz f(\ba)$.  However,
$f_\fz(\ba)$ is also well-defined on any $\ba\in\tlgA$ such
that $\ba_{\dB(\fz,r)}$ is $\Sft$-admissible.  
Hence $f_\fz(\ba)$ is well-defined for all $\fz\in\unflawed(\ba)$.
If $\ba\in\tlgA$ and $\fz,\fy\in\unflawed(\ba)$,
then we say that $\ba$ has a {\dfn $\blam$-dislocation
between $\fy$ and $\fz$} if
$f_\fy(\ba) \neq \blam^{\fy-\fz} f_\fz(\ba)$.  Let
\beqn
\label{spec.gap}
  \gam_{\fy,\fz}(\blam)
\quad:=\quad 
\blam^{\fy-\fz} \cdot f_\fz(\ba) /  f_\fy(\ba).  
\eeqn
If $\SFT$ is transitive, then $\gam_{\fy,\fz}(\blam)$ is
independent of the choice $f\in\Eigenspace$, by Lemma 
\ref{spectral.lemma}(b).

\Theorem{\label{gap.prop}}
{
Let $\SFT$ be transitive.   Let $\blam\in\RatSpec$
{\rm (with radius $r'>0$)} and 
let $\ba\in\tlgA$ have a $\blam$-dislocation.  Let $r\geq r'+1$.  Then:
\bthmlist 
  \item $\ba$ has a range $r$ domain boundary.

  \item Suppose $\unflawed(\ba)$ has connected components
 $\{\dY_n\}_{n=1}^N$, where $N\in\Natur\union\{\oo\}$.
There is a matrix of displacements
 $\Del:=[\del_{nm}]_{n,m=1}^N$ such that
 \bitem
  \item[{\rm[i]}] $\gam_{\fy,\fz}(\blam) \ = \ \del_{nm}(\blam)$
for any $\fy\in\dY_n$ and $\fz\in\dY_m$.

  \item[{\rm[ii]}] {\rm(Cocycle property)} \ $\del_{n \ell}(\blam) = \del_{nm}(\blam)\del_{m \ell}(\blam)$ for any  $n,m,\ell\in\CC{1...N}$.
\eitem

  \item If two of the components $\dY_n$ and $\dY_m$ are projective,
and $\del_{nm}$ is nontrivial, then there is an
essential domain boundary between $\dY_n$ and $\dY_m$
\ethmlist
}
\bthmprf {\bf(a)} follows from {\bf(b)}, because if $\unflawed(\ba)$
had only one
connected component $\dY_1$, then {\bf(b)}[ii] implies that
$\del_{11}\equiv 1$, and then  {\bf(b)}[i] says
that $\gam_{\fy,\fz}(\blam)  =  1$
for all $\fy,\fz\in \dY_1$, which contradicts the hypothesis that $\ba$
has a $\blam$-dislocation.  To prove {\bf(b)} we need the following:

 \Claim{\label{gap.lemma}
{\rm[a]} \ $\gam_{\fx,\fz}(\blam)=1$
if $\fx,\fz$ are in the same connected component
of $\unflawed(\ba)$.

{\rm[b]} \ For any $\fx,\fy,\fz\in \unflawed(\ba)$, we have
$\gam_{\fx,\fz}(\blam) \ = \ 
\gam_{\fx,\fy}(\blam)\cdot \gam_{\fy,\fz}(\blam)$.

{\rm[c]} \ If $\blam_1,\blam_2\in\RatSpec$ have radius $r$,
then $\blam:=\blam_1\cdot\blam_2$ also has radius $r$, and
$\gam_{\fx,\fy}(\blam) \ = \ \gam_{\fx,\fy}(\blam_1)
\cdot \gam_{\fx,\fy}(\blam_2)$.
}
\bclaimprf
{\rm[a]} Let $\fx,\fz\in\dY_n$ for some $n\in\CC{1...N}$.
Let $\fx=\fy_0\adjacent\fy_1\adjacent\cdots\adjacent\fy_J=\fz$ be
a trail in $\unflawed$ from $\fx$ to $\fz$.
Then for all $j\in\CC{1...J}$, \ $\kap_{\fy_j}(\ba)=\blam^{\fy_j-\fy_{j-1}}
\kap_{\fy_{j-1}}(\ba)$
({\em Proof:} Similar to Claim 1 of Proposition \ref{interface.prop}).
Inductively, we have
$f_{\fz}(\ba) \ = \ \blam^{\fz-\fx} f_{\fx}(\ba)$.  Thus,
 $\gam_{\fx,\fz}(\blam)  =  1$.

{\rm[b]} follows from eqn.(\ref{spec.gap}). For
{\rm[c]}, let $f_1\in\Eigenspace[\blam_1]$ and 
 $f_2\in\Eigenspace[\blam_2]$, and let $f:=f_1\cdot f_2$.
Lemma \ref{spectral.lemma}(a) says that 
$f\in\Eigenspace$; hence radius$(\blam)=r$.
Now substitute $f_1$, $f_2$, and $f$ respectively
into eqn.(\ref{spec.gap}) to compute $\gam_{\fy,\fz}(\blam_1)$,
\ $\gam_{\fy,\fz}(\blam_2)$ and $\gam_{\fy,\fz}(\blam)$.
\eclaimprf
{\bf(b)}[i] follows Claim 1[a,b].
Then {\bf(b)}[ii] is by  Claim 1[b].
Finally, Claim 1[c] implies that \linebreak $\del_{n,m}:\RatSpec\into\dT$ is
a homomorphism (i.e. a displacement).

{\bf(c)}  If $\dY_1$ and $\dY_2$ are projective,
then we cannot eliminate
the dislocation by changing $\ba$ inside $\unflawed[R](\ba)$ for
any $R\in\Natur$, so the defect is essential.\ethmprf

Theorem \ref{gap.prop}(b) allows us to speak of a {\dfn rational $(\gA,\Phi)$-dislocation}  rather than of a ``$\blam$-dislocation''.
(When $\gA$ and $\Phi$ are clear, we will just say ``rational dislocation'').
The cocycle property means: {\bf(a)} All entries of the
displacement matrix
$\Del_{\ba}$ can be reconstructed from one row, and
{\bf(b)}  $\Del_{\ba}$
is `antisymmetric', i.e. $\del_{nm}(\blam) = \del_{mn}(\blam)^{-1}$
and $\del_{nn} \equiv 1$, $\forall \ n,m\in\CC{1...N}$.
If there are only two connected
components, $\dY_1$ and $\dY_2$, then we define the
{\dfn displacement} of $\ba$ to be $\del_\ba:=\del_{12}$.

\newcommand{\domain}[1]{\fbox{$#1$}}

\begin{figure}[h]
\centerline{
\begin{tabular}{|c|c|c|c|}
\hline
\includegraphics[height=5.5em,width=6.5em]{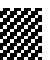} &
\includegraphics[height=5.5em,width=6.5em]{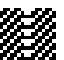} &
\includegraphics[height=5.5em,width=6.5em]{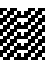} &
\includegraphics[height=5.5em,width=6.5em]{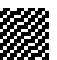} \\
\hline
$(*)$ & $\boldsymbol{(\alpha)}$ & $\boldsymbol{(\beta)}$  &  $\boldsymbol{(\gamma)}$ \\
\hline
\end{tabular} }
\Caption{fig:62.defects}{$(*)$  The periodic background generated by $\ECA{62}$
acting on $\gD$;\quad
$\boldsymbol{(\alpha,\beta,\gamma)}$: Three rational dislocations;
see Examples \ref{X:persistent.dislocation}(a)
and \ref{X:permanent.dislocations}(a).
(In \cite[\S III(C)]{BNR91}, $\alp$, $\beta$, and $\gam$ are 
``$g_e$'', ``$g_o$'' and ``$w$'' respectively.)}
\end{figure}

\EXAMPLE{\label{X:persistent.dislocation}
  (a) (ECA\#62) Let $\sA=\{\black,\white\}$, and let $\gD$ be the
three-element $\shift{}$-orbit of $\overline{\black\black\white}$.
Let $\lam:=e^{2\pi\bi/3}$.
Then $\gD$ is $\shift{}$-transitive, and 
 $\RatSpec[\gD,\shift{}]=\{1,\lam,\lam^2\}\subset\dT$
[because $\gD$ has period $P=3$; see Example \ref{X:ratspec}(a)].
If $\ECA{62}$ is ECA\#62, then $\ECA{62}\restr{\gD}=\shift{}$
[see Figure \ref{fig:62.defects}$(*)$], so
the homomorphism $\tau:\Spec[\gD,\shift{}]\into\dT$ in 
Lemma \ref{spectral.lemma}{(d)} is the identity:
$\tau(\lam^p)=\lam^{p}$ [see Example \ref{X:ratspec}(a)].
Thus, 
\[
\RatSpec[\gD,\ECA{62},\shift{}]\quad=\quad\{(1,1),(\lam,\lam),(\lam^2,\lam^2)\}
\quad\subset\quad\dT^2.
\]
  The homomorphism $\Zahl\ni\fz\mapsto\displace{(0;\fz)}\in
\hRatSpec[\gD,\ECA{62},\shift{}]\cong\Zahlmod{3}$ is surjective,
with kernel $3\Zahl$.  Hence we identify displacements with
elements of $\Zahlmod{3}$.
Below are three rational dislocations in $\gD$ and their displacements.

\centerline{$
\begin{array}{rlr}
\alp & \domain{\black\black\white  \black\black\white}  
\overrightarrow{ \ \white\white\white \ }  \domain{\black\black\white  \black\black\white} & 
\del=3\equiv 0\in\Zahlmod{3} \\
\bet & \domain{\black\black\white  \black\black\white}  
\overrightarrow{ \ \white\white \ }   \domain{\black\black\white   \black\black\white\black} & 
\del=2\in\Zahlmod{3} \\
\gam & \domain{\black\black\white   \black\black\white}  
\overrightarrow{\ \white \ }  \domain{\black\black\white  \black\black\white\black\black} & 
\del=1\in\Zahlmod{3} 
\end{array}$}

[compare to the top rows in  Figure \ref{fig:62.defects}$(\alp,\bet,\gam)$.]
The $\gam$ and $\bet$ dislocations have nontrivial
displacements, so they are are essential.
The $\alp$ dislocation is not essential (it can
be removed by replacing the middle three blocks with $\black\black\white$).

\begin{figure}[h]
\centerline{
\begin{tabular}{|c|c|c|c|c|}
\hline
\includegraphics[height=5.5em,width=6em]{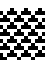} &
\includegraphics[height=5.5em,width=6em]{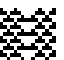} &
\includegraphics[height=5.5em,width=6em]{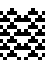} &
\includegraphics[height=5.5em,width=6em]{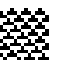} &
\includegraphics[height=5.5em,width=6em]{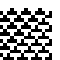} \\
\hline
$(*)$  & $\boldsymbol{(\alpha)}$ & $\boldsymbol{(\beta)}$ &  $\boldsymbol{(\gamma^+)}$ & $\boldsymbol{(\gamma^-)}$  \\
\hline 
\end{tabular} }
\Caption{fig:54.defects}{$(*)$  The periodic background generated by $\ECA{54}$
acting on $\gB$.\quad
$\boldsymbol{(\alpha, \beta,\gamma^\pm)}$: Four rational
dislocations;  see Examples \ref{X:persistent.dislocation}(b)
and \ref{X:permanent.dislocations}(b).
(This nomenclature is due to \cite[Fig.8]{CrHa97}.
In the nomenclature of \cite[\S III(C)]{BNR91}, 
$\alp=g_o$, $\bet=g_e$, \ 
$\gamma^+=\vec{w}$ and $\gamma^-=\stackrel{\leftarrow}{w}$.)
See also \cite[Example 1.2(b)]{PivatoDefect0}.
}
\end{figure}

 (b) (ECA\#54)  Let $\sA=\{\black,\white\}$, and let $\gB:=\gB_0\disj\gB_1$,
where $\gB_0$ is the four-element $\shift{}$-orbit of $\overline{\black\black\black\white}$ and 
$\gB_1$ is the  four-element $\shift{}$-orbit of $\overline{\white\white\white\black}$.  
Then $\gB$ is not $\shift{}$-transitive,
and \[
\RatSpec[\gB_0,\shift{}]
\quad=\quad\RatSpec[\gB_1,\shift{}]
\quad=\quad\{\bi^p\}_{p=0}^3
\quad=\quad\{1,\bi,-1,-\bi\}
\quad\subset\quad\dT
\]
 [because $\gB$ has period $P=4$;
see Example \ref{X:ratspec}(a)].  
If $\ECA{54}$ is ECA\#54, then $\ECA{54}(\gB_0)=\gB_1$ and
$\ECA{54}(\gB_1)=\gB_0$, so $\gB$ is $(\ECA{54},\shift{})$-transitive.
Also, $\ECA{54}^2\restr{\gB}=\shift{2}$ 
[see Figure \ref{fig:54.defects}$(*)$],
so the epimorphism $\Zahl\x\Zahl\ni(t;\fz)\mapsto\displace{(t;\fz)}\in
\hRatSpec[\gB,\ECA{54},\shift{}]$
has kernel $\dK:=\Zahl(2,2) \dirsum \Zahl(0,4)$.  Hence we identify displacements with
elements of $\ZD[2]/\dK$.
Below are four  rational dislocations and their displacements
 [compare to the top rows in Fig.\ref{fig:54.defects}$(\alp,\bet)$, 
or middle rows in Fig.\ref{fig:54.defects}$(\gam^\pm)$].

\centerline{\footnotesize
\psfrag{a0}[][]{$\alp_0$}
\psfrag{a1}[][]{$\alp_1$}
\psfrag{a2}[][]{$\alp_2$}
\psfrag{a3}[][]{$\alp_3$}
\psfrag{b0}[][]{$\bet_0$}
\psfrag{b1}[][]{$\bet_1$}
\psfrag{b2}[][]{$\bet_2$}
\psfrag{b3}[][]{$\bet_3$}
\psfrag{sa}[][]{$\del=(0,3)+\dK$}
\psfrag{sb}[][]{$\del=(0,2)+\dK$}
\includegraphics[scale=0.48,angle=-90]{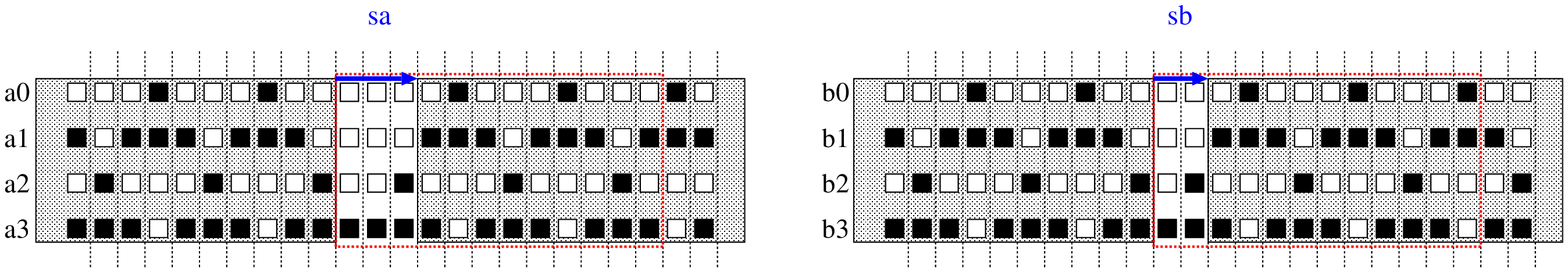}
}

\centerline{\footnotesize
\psfrag{g0}[][]{$\gam_0^+$}
\psfrag{g1}[][]{$\gam_1^+$}
\psfrag{G0}[][]{$\gam_0^-$}
\psfrag{G1}[][]{$\gam_1^-$}
\psfrag{sgp}[][]{$\del=(1,1)+\dK$}
\psfrag{sgm}[][]{$\del=(-1,1)+\dK$}
\includegraphics[scale=0.48,angle=-90]{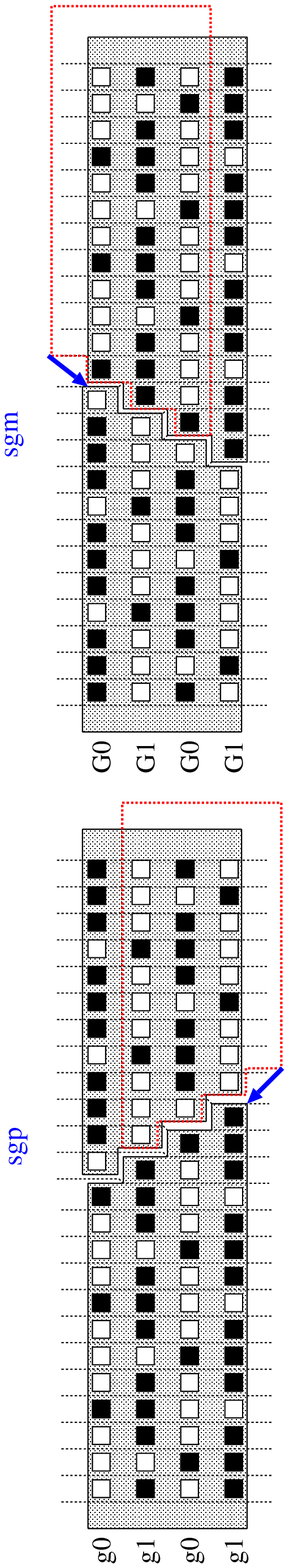}
}

\vspace{0.5em}

(c) (ECA\#184) Let $\sA=\{\black,\white\}$, and let $\gG_*\subset\AZ$ be 
as in Example  \ref{X:persistent.interface}(a). Then
$\RatSpec[\gG_*,\shift{}]=\{\pm1\}$ 
[because $\gG_*$ has period 2; see Example \ref{X:ratspec}(a)].  
Also, 
$\ECA{184}\restr{\gG_*}=\shift{}$ [see Figure \ref{fig:184.defects}$(*)$], so
the homomorphism \linebreak $\tau:\Spec[\gG_*,\shift{}]\into\dT$ in 
Lemma \ref{spectral.lemma}{(d)}  is the identity
[see Example \ref{X:ratspec}(a)]. 
Thus, $\Spec[\gB_*,\ECA{184},\shift{}]=\{(1,1),(-1,-1)\}$,
and the homomorphism $\Zahl\ni\fz\mapsto\displace{(0;\fz)}\in
\hRatSpec[\gG_*,\ECA{184},\shift{}]\cong\Zahlmod{2}$ is surjective,
with kernel $2\Zahl$, so we identify displacements with
elements of $\Zahlmod{2}$.
Below are two rational dislocations and their displacements
[compare to any rows in Figure \ref{fig:184.defects}$(\gam^\pm)$]:
\[
\begin{array}{rll}
\gam^+ & \domain{\white\black\white\black\white\black\white\black} 
\overrightarrow{ \ \black \ }
\domain{\white\black\white\black\white\black\white\black}
  \  &  \del=1\in\Zahlmod{2} \\
\gam^- & 
\domain{\white\black\white\black\white\black\white\black}
\overrightarrow{ \  \white \ }
\domain{\white\black\white\black\white\black\white\black} \
 &  \del=1\in\Zahlmod{2}
\end{array}
\]
\begin{sloppypar}
(d) (ECA \#110) Let $\sA:=\{\black,\white\}$, and let $\gE\subset\AZ$ be
the 14-element $\shift{}$-orbit of the 14-periodic sequence
$\overline{\black\black\black\white\black\black\white\white\black\white\white\white\white\white}$.  
Let $\lam:=e^{\pi\bi/7}$.
Then $\RatSpec[\gE,\shift{}]=\{1,\lam,\ldots,\lam^{13}\}\subset\dT$
[because  $\gE$ has period $14$; see Example \ref{X:ratspec}(a)].
If $\ECA{110}$ is ECA\#110, then $\ECA{110}\restr{\gE}=\shift{4}$
[see Figure \ref{fig:110.defects}$(*)$], so
the homomorphism $\tau:\Spec[\gE,\shift{}]\into\dT$ in 
Lemma \ref{spectral.lemma}{(d)}  is given
$\tau(\lam^p)=\lam^{4p}$ [see Example \ref{X:ratspec}(a)].
Thus,
 \[
\Spec[\gD,\ECA{62},\shift{}]\quad=\quad\{(1,1),(\lam^4,\lam),(\lam^{8},\lam^2),
\ldots,(\lam^{10},\lam^{13})\},
\]
and the homomorphism $\Zahl\ni\fz\mapsto\displace{(0;\fz)}\in
\hRatSpec[\gE,\ECA{110},\shift{}]\cong\Zahlmod{14}$ is surjective,
with kernel $14\Zahl$, so we identify displacements with
elements of $\Zahlmod{14}$.
Figure \ref{fig:110.defectsb} shows
seven essential rational dislocations in $\gE$ with
nontrivial displacements. 
Figure \ref{fig:110.defects} shows their $\ECA{110}$-evolution.
\end{sloppypar}

\begin{figure}[h]
{\scriptsize \psfrag{sA}[][]{$\del=6\in\Zahlmod{14}$}
\psfrag{sB}[][]{$\del={8}\in\Zahlmod{14}$}
\psfrag{sC}[][]{$\del=9\in\Zahlmod{14}$}
\psfrag{sD1}[][]{$\del={11}\in\Zahlmod{14}$}
\psfrag{sE}[][]{$\del=23\equiv{9}\in\Zahlmod{14}$}
\psfrag{sEb}[][]{$\del={5}\in\Zahlmod{14}$}
\psfrag{sF}[][]{$\del=15\equiv 1\in\Zahlmod{14}$}
\includegraphics[width=12em,height=50em,angle=-90]{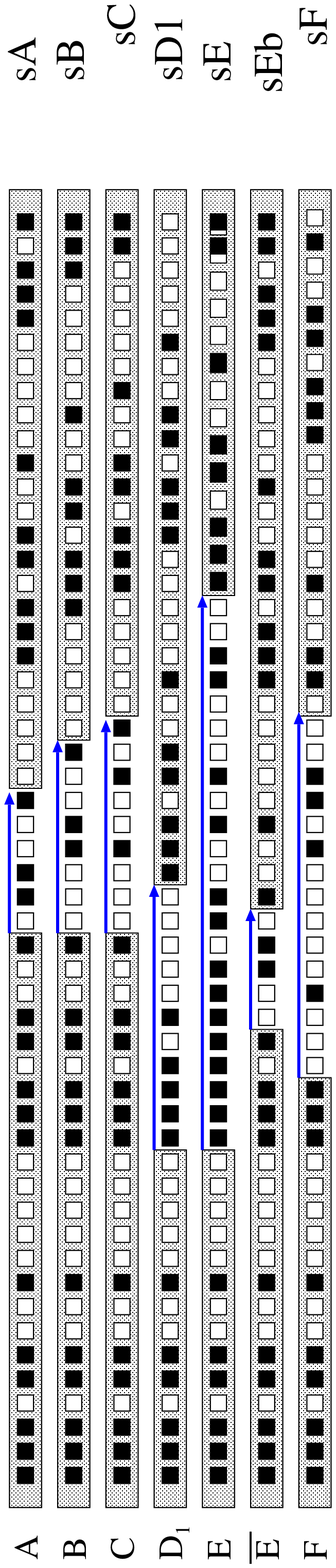}
}
\Caption{fig:110.defectsb}{Seven essential rational dislocations in $\gE$ with
nontrivial displacements; see Examples \ref{X:persistent.dislocation}(d)
and \ref{X:permanent.dislocations}(d).} 
\end{figure}

\ignore{\[
\psfrag{sA}[][]{$\del=6\in\Zahlmod{14}$}
\psfrag{sB}[][]{$\del={8}\in\Zahlmod{14}$}
\psfrag{sC}[][]{$\del=23\equiv 9\in\Zahlmod{14}$}
\psfrag{sD1}[][]{$\del={11}\in\Zahlmod{14}$}
\psfrag{sE}[][]{$\del={9}\in\Zahlmod{14}$}
\psfrag{sEb}[][]{$\del={5}\in\Zahlmod{14}$}
\psfrag{sF}[][]{$\del=15\equiv 1\in\Zahlmod{14}$}
\includegraphics[width=8em,height=38em,angle=-90]{Pictures/110/dislocate.eps}
\]}

(e) \label{X:dislocation}
 Let $\sA=\{\black,\white\}$ and let $\Checker\subset\AZD[2]$ be
as in Example \ref{X:codim.one}(b).
Then $\RatSpec[\Checker,\shift{}] = \{\pm 1\}^2\subset\dT^2$,
by Example  \ref{X:ratspec}(b) [because $\Checker$ is $\dP$-periodic, where
$\dP:=\Zahl(1,1)\dirsum\Zahl(1,-1)$].
Let $\Phi$ be a CA with $\Checker\subseteq\Fix{\Phi}$
[e.g. Example \ref{X:biject.persistent.defect}(b)].
Then $\RatSpec[\Checker,\Phi,\shift{}]$ $ = \{(1;\pm 1,\pm1)\}$ $\subset
\dT^3$ and $\hRatSpec[\Checker,\Phi,\shift{}]\cong\ZD[2]/\dP$.
 [by Lemma \ref{spectral.lemma}(d) and Example \ref{X:ratspec}(b)].
 Figure \ref{fig:codim}(B) shows a domain boundary in $\Checker$.  If
$\dY_0$ is the northern connected component and $\dY_1$ is the
southern component, then we have $\del_{12}=(1,0)+\dP$.  Both
components are projective, so this is an essential rational dislocation.}

\begin{figure}[h]
\centerline{
\begin{tabular}{|c|c|c|c|c|c|c|c|}
\hline
\includegraphics[height=5.5em,width=7em]{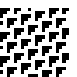} &
\includegraphics[height=5.5em,width=7em]{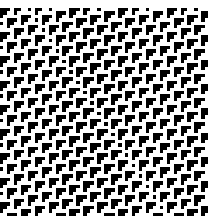} &
\includegraphics[height=5.5em,width=7em]{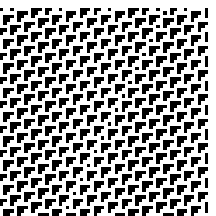} &
\includegraphics[height=5.5em,width=7em]{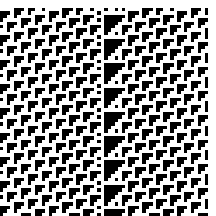} \\
$(*)$     & {\bf(A)}   & {\bf(B)} & {\bf(C)} \\
\includegraphics[height=5.5em,width=7em]{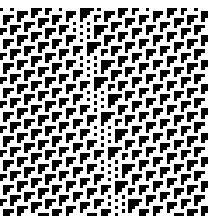} &
\includegraphics[height=5.5em,width=7em]{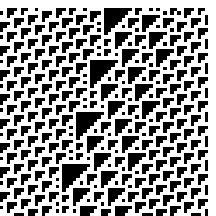} &
\includegraphics[height=5.5em,width=7em]{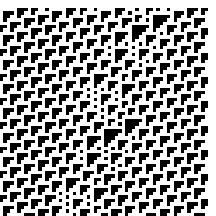} &
\includegraphics[height=5.5em,width=7em]{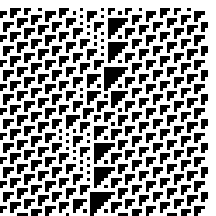} \\
$\mathbf{(D_1)}$ & {\bf(E)} {\scriptsize(`extended')}  & $\mathbf{(\barE)}$  &  {\bf(F)} \\
\hline
\end{tabular} }
\Caption{fig:110.defects}{$(*)$  A $30\x30$ image of
the periodic spacetime diagram of $\ECA{110}$ acting on $\gE$;\quad
{\bf(A,B,C,$\mathbf{D_1}$,E,$\mathbf{\barE}$,F):}  $60\x60$ images of the
$\ECA{110}$-evolution of seven dislocations in $\gE$. See Examples \ref{X:persistent.dislocation}(d) and \ref{X:permanent.dislocations}(d).
(This nomenclature is due to \cite{Cook}.
In  the nomenclature of \cite{CrHR} $A=\omg_{^{\mathrm{right}}}$, $B=\omg_{^{\mathrm{left}}}$, $C=\alp$, $E=\bet$, $F=\eta$ etc.).
See also \cite[Example 1.2(c)]{PivatoDefect0}.}
\end{figure}

\Proposition{\label{one.dim.SFT.essential.persistent}}
{
  Let $\gA\subset\AZ$ be a one-dimensional SFT.
Let $\ba\in\tlgA$.
\bthmlist
\item If $\gA$ is $\shift{}$-mixing, then every finite defect of $\ba$ is removable.

\item If $\gA$ is $(\Phi,\shift{})$-transitive (but not $\shift{}$-mixing), 
then any finite essential defect in $\ba$ is a rational dislocation.

\item If $\gA$ is not $(\Phi,\shift{})$-transitive, then any finite  essential
defect is either a rational dislocation or an interface.
\ethmlist
}
\bthmprf  {\bf(a)} is Example \ref{X:removable.finite.defects}(a). 
{\bf(c)} follows from {\bf(b)}, which
follows from Example \ref{X:one.dim.per.SFT}.
\ethmprf

%\subsection{Persistence and Coalescence of Rational Dislocations}

%  Now we can address questions \#\ref{Q:persist} and \#\ref{Q:chem} from the introduction.
Next we show that any essential, rational $(\gA,\Phi$)-dislocation is
$\Phi$-persistent:

\Theorem{\label{persistent.dislocation}}
{
Let $\Phi:\AZD\into\AZD$ be a CA and let $\gA\subset\AZD$ be a subshift
with $\Phi(\Sft)=\Sft$. 
If $\ba\in\tlgA$ has an essential
rational dislocation, then
$\Phi(\ba)$ also has an essential rational dislocation,
with the same displacement matrix  as $\ba$.
}
\bthmprf 
Let $\ba':=\Phi(\ba)$, and suppose $\Phi$ has radius $R>0$.
Then each projective component of $\unflawed[r+R](\ba)$ is
contained in a projective component of $\unflawed(\ba')$,
because $\unflawed[r+R](\ba)\subseteq\unflawed[r](\ba')$
by Lemma \ref{defect.dimension}(b).
 Let $\unflawed[R+r](\ba)$  have projective components
$\dY_1,\ldots,\dY_N$, and let
$\unflawed(\ba')$  have projective components
$\dY'_1,\ldots,\dY'_N$, where $\dY_n\subset\dY'_n$ for all $n\in\CC{1...N}$.
Let $\Del_{\ba} = [\del_{nm}]$ and
$\Del_{\ba'}=[\del'_{nm}]$.
We must show $\del_{nm}=\del'_{nm}$ for all $n,m\in\CC{1...N}$.
If $\lam\in\RatSpec$
and  $f\in\Eigenspace$, then $f\circ\Phi\,\eeequals{ae}\,\lam_0 f$.
Let $\fy\in\dY_n\subset\dY'_n$ and let 
$\fz\in\dY_m\subset\dY'_m$.  Then
\[
  \del_{nm}(\blam) \quad\eeequals{(*)}\quad
\blam^{\fy-\fz} \cdot \frac{f_\fz(\ba)}{ f_\fy(\ba)}
\quad=\quad
\blam^{\fy-\fz} \cdot \frac{\lam_0 f_\fz(\ba)}{\lam_0 f_\fy(\ba)}
\quad\eeequals{(\dagger)}\quad
\blam^{\fy-\fz} \cdot \frac{f_\fz(\ba')}{f_\fy(\ba')}
\quad\eeequals{(*)}\quad
\del'_{nm}(\blam).
\]
Here, $(*)$ is by eqn.{\rm(\ref{spec.gap})} and 
Theorem \ref{gap.prop}(b)[i], and
$(\dagger)$ is because $\lam_0 f(\ba)  =  (f\circ\Phi)_\fz(\ba)
 =  f_\fz(\ba')$. 
\ethmprf

\EXAMPLE{\label{X:permanent.dislocations}
(a) The $\gam$ and $\bet$ dislocations of Example \ref{X:persistent.dislocation}(a) are $\ECA{62}$-persistent, by Theorem \ref{persistent.dislocation}.
(The $\alp$ dislocation is also $\ECA{62}$-persistent, but not because of 
Theorem \ref{persistent.dislocation}). The $\ECA{62}$-evolution of
all three dislocations is shown in 
Figure \ref{fig:62.defects}$(\alp,\bet,\gam)$.
Their large-scale $\ECA{62}$-dynamics were shown in Figure \ref{fig:defect.intro}(C).

(b) All four dislocations in Example \ref{X:persistent.dislocation}(b)
are $\ECA{54}$-persistent; their $\ECA{54}$-evolution is shown in 
 Figure \ref{fig:54.defects}$(\alp,\bet,\gam^\pm)$.  See also  Figure \ref{fig:defect.intro}(B).

(c) Both dislocations in Example \ref{X:persistent.dislocation}(c)
are $\ECA{184}$-persistent; 
their $\ECA{184}$-evolution is shown in Figure
\ref{fig:184.defects}$(\gam^\pm)$.
See also  Figure \ref{fig:defect.intro}(D).

(d) All seven dislocations  in Example
\ref{X:persistent.dislocation}(d) are $\ECA{110}$-persistent.  Figure
\ref{fig:110.defects} shows their $\ECA{110}$-evolution.  These are
only some of the plethora of defect particles of ECA \#110 (see
\cite{Lind2,McIntosh1,CrHR} or \cite[\S3.1.4.4]{Ilachinski}), whose
complex interactions can support universal computation;
see \cite{Cook}, \cite{McIntosh2} or \cite[Ch.11]{Wolfram2}.

(e)  If  $\gA\subset\AZ$ is any one-dimensional SFT,
then every essential defect is persistent.  (Combine
Propositions \ref{persistent.interface} and
\ref{one.dim.SFT.essential.persistent}(c) with
Theorem \ref{persistent.dislocation}). 
\projdisl{

(f)  The essential dislocation in Example \ref{X:irrational.dislocation}(a)
is $\ECA{18}$-persistent.}
}

\begin{figure}[h]
\centerline{\scriptsize
\begin{tabular}{|c|c|c|c|c|}
\hline
\multicolumn{2}{|c|}{{\normalsize ECA \#62}} &{\normalsize ECA \#184} & \multicolumn{2}{|c|}{{\normalsize ECA \#54}}\\
\hline
\includegraphics[width=7.8em,height=8em]{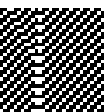} &
\includegraphics[width=7.8em,height=8em]{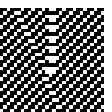} &
\includegraphics[width=7.8em,height=8em]{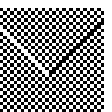} &
\includegraphics[width=7.8em,height=8em]{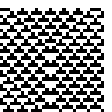} &
\includegraphics[width=7.8em,height=8em]{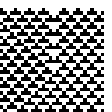} \\
{\bf(a)} $\gam + \bet \rightarrow \alp$ &
{\bf(b)}  $\gam + \alp \rightarrow \gam$ &
{\bf(c)}  $\gam^+ + \gam^- \rightarrow \emptyset$ &
{\bf(d)}  $\gam^+ + \gam^- \rightarrow \beta$ &
{\bf(e)}  $\gam^+ + \bet \rightarrow \gam^-$ \\
$2 + 1 \equiv 0\in\Zahlmod{3}$&
$2 + 0 \equiv 2\in\Zahlmod{3}$&
$1+1 \equiv 0\in\Zahlmod{2}$ &
$(1,1)+(-1,1)=$&
$(1,1)+(0,2)=$\\
 & & & $(0,2)\in\ZD[2]/\dK$&
$(1,3)\equiv(-1,1)\in\ZD[2]/\dK$\\
\hline
\end{tabular} }
\Caption{fig:collisions}{Dislocation coalescence and the algebra of $\hSpecae$.}
\end{figure}

\parag{Defect coalescence:}
If $D=1$, then dislocations can be thought of as `particles';
see \cite{PivatoDefect0}.
If two such `dislocation particles' $x$ and $y$ coalesce to form $z$, then
$\del_z =\del_x \cdot \del_y$.  In particular, $x$ and $y$ can
annihilate only if $\del_x \cdot \del_y\equiv 1$.
Thus, the algebra of $\hSpecae$ yields `conservation laws' which
helps to determine the `chemistry' of dislocation particles.
This partially answers Question \#\ref{Q:chem} from the introduction.
See Figure \ref{fig:collisions} for some examples.

%spellchecked
%\newcommand{\brv}[1]{\breve{#1}}
\newcommand{\undgA}{\underline{\gA}}
\newcommand{\undgS}{\underline{\gS}}
\newcommand{\undgT}{\underline{\gT}}

\subsection{\label{S:projective.dislocation}Projective Dislocations}

If $\gA\subset\AZ$ is a one-dimensional subshift of finite type, then Proposition
\ref{one.dim.SFT.essential.persistent} completely characterizes its
essential defects.  However, if $\gA$ is {\em not} of finite type,
or even if $\gA\subset\AZD$ {\em is} of finite type, but $D\geq 2$, then Proposition
\ref{one.dim.SFT.essential.persistent} fails, because some
$\shift{}$-eigenvalues of $\gA$ may be irrational\footnote{This is can
occur in higher dimensional SFTs which are conjugate to substitution
shifts \cite{Mozes}; see
\cite[\S3(d)]{Radin} or \cite[\S7.1]{Robinson} for a discussion.}, and
even rational eigenvalues may only have discontinuous eigenfunctions.
Theorem \ref{gap.prop} is not applicable to such eigenfunctions,
because they are {\em not} locally determined.  We can extend
the methods of \S\ref{S:dislocation} 
to irrational or \ae-continuous eigenfunctions, but only for
projective domain boundaries, and only if $\gA$ satisfies certain
transitivity and extendibility conditions. 

  A {\dfn meager} subset of $\gA$ is any countable union of closed,
nowhere-dense sets.  A {\dfn comeager} (or {\dfn residual}) subset of
$\gA$ is the compliment of a meager subset; the family of {comeager}
subsets of $\gA$ is thus closed under countable
intersections and arbitrary unions.  The Baire Category Theorem
\cite[Thm.5.8]{Folland} says that any comeager set is dense in $\gA$.
A bounded function $f\colon\gA\rightarrow\Cplx$ is {\dfn almost
everywhere {\rm(``\ae'')} continuous} if $f$ is continuous at each point in a
comeager subset of $\gA$.  If $f,g\colon\gA\rightarrow\Cplx$, then we
say that $f = g$ {\dfn almost everywhere} (``$f \ \eeequals{ae}\ g$'')
if $\set{\ba\in\gA}{f(\ba)=g(\ba)}$ is comeager in $\gA$.  This is an
equivalence relation. Let $\Cae(\gA)$ be the $\Cplx$-vector space of
\ae-equivalence classes of \ae-continuous functions from $\gA$ into
$\Cplx$.  An element $f\in\Cae(\gA)$ is an (\ae) {\dfn
$(\Phi,\shift{})$-eigenfunction} with (\ae) {\dfn eigenvalue} $\blam =
(\lam_0; \lam_1, ...,\lam_D)\in\Torus{D+1}$ if: 
\bdesc
\item[{\bf(a)}]
$f\circ\Phi \ \eeequals{ae} \ \lam_0 f$.
\item[{\bf(b)}] For any
$\fz=(z_1,\ldots,z_D)\in\ZD$, $f\circ \shift{\fz} \ \eeequals{ae} \
\blam^\fz f$.
\edesc Let $\Specae\subset\Torus{D+1}$ be the set of all such
eigenvalues.  For any $\blam\in\Specae$, let 
\[
\Eigenspace \quad:=\quad
\set{f\in\Cae(\Sft)}{f\circ\Phi \ \eeequals{ae} \  \lam_0 f \ \AND \  f\circ\shift{\fz} \ \eeequals{ae} \ 
\blam^\fz f,
\ \forall  \fz\in\ZD}
\] be the {\dfn eigenspace} of
$\blam$.

\EXAMPLE{\label{X:ae.continuous}
(a) Let $\gS$ be the sofic shift of Example \ref{X:defect.energy}(c).
Let \[
\gE \ := \ \{\ba\in\AZ \ ; \ a_z = 0, \ \forall \ \mbox{even} \
z\in\Zahl\}\quad\AND\quad \gO \ := \ \{\ba\in\AZ \ ; \ a_z = 0, \ 
\forall \ \mbox{odd}  \ z\in\Zahl\}.
\]
Then $\gS:=\gE\union\gO$, and $\gE\intsct\gO=\{\bar0\}$.
Let $\gE^*:=\gE\setminus\{\bar0\}$ and $\gO^*:=\gO\setminus\{\bar0\}$.
Then $\shift{}(\gE^*)=\gO^*$ and $\shift{}(\gO^*)=\gE^*$.
Define $f:\gS\into\{-1,0,1\}$ by $f\restr{\gO^*}\equiv -1$, \ 
$f(\bar0)= 0$, \  and $f\restr{\gE^*}\equiv 1$.  Then $f\in\Cae(\gS)$
and is a $\shift{}$-eigenfunction with eigenvalue $-1$. 
However, $\Eigenspace[-1](\gS,\shift{})$ has no {\em continuous}
eigenfunctions, because $\gS$ contains a dense subset of
 points which are $\shift{}$-homoclinic to $\bar0$ (namely sequences with only
finitely many ones).  
Hence $-1\in\Specae[\gS,\shift{}]\setminus \Spec[\gS,\shift{}]$.

  If $\ECA{18}$ is ECA\#18, then $\ECA{18}(\gE)=\gO$ and 
$\ECA{18}(\gO)=\gE$. Hence $f\in \Eigenspace[(-1;-1)](\gS,\Phi,\shift{})$,
and $(-1;-1)\in\Specae[\gS,\ECA{18},\shift{}]$.  In this case,
\[
\Specae[\gS,\ECA{18},\shift{}]\quad=\quad\{(1,1),(-1,-1)\}\quad\subset\quad
\dT^2,
\]
and the function $\Zahl \ni \fz\mapsto\displace{(0;\fz)}\in
\hSpecae[\gS,\ECA{18},\shift{}]\cong\Zahlmod{2}$ is a surjection
with kernel $2\Zahl$; hence we identify displacements with 
elements of $\Zahlmod{2}$.

(b) Let $\sA=\{0,1\}$, let $\lam\in\CO{0,1}$ be  irrational, and
define $\tau:\CO{0,1}\ni x\mapsto(x+\lam \bmod 1)\in\CO{0,1}$
(i.e. the rotation system induced by $\lam$).
Let $\gT\subset\AZ$ be the Sturmian subshift \cite[Ch.6]{Fogg}
obtained by projecting $\tau$-orbits through the partition
$\sP:=\{\CO{0,1-\lam}, \ \CO{1-\lam,1}\}$.
  Then $\Specae[\gA,\shift{}]=
\Spec[\gA,\shift{}]=\{e^{2n\pi\bi\lam}\}_{n\in\Zahl}\subset\dT$,
so $\RatSpec[\gA,\shift{}]$ is trivial.  In this case,
the function $\Zahl\ni \fz\mapsto \del_{\vec{\fz}}\in\hSpecae[\gA,\shift{}]$
is an isomorphism, so we can identify all displacements with integers.
}
The system $(\gA,\Phi,\shift{})$ is 
{\dfn topologically weakly mixing} if the Cartesian product 
$(\gA\x\gA, \Phi\x\Phi,\shift{}\x\shift{})$ is topologically transitive.

\Lemma{\label{keynes.robertson}}
{
{\rm[i]} $\Spec\subseteq\Specae$.  If $(\gA,\Phi,\shift{})$ is transitive, then
{\rm Lemma \ref{spectral.lemma}{\bf(a-e)} [but not {\bf(f)}]} are still true if
we replace ``$\sC(\gA)$'' with ``$\Cae(\gA)$'' and 
``$\Spec$'' with ``$\Specae$''. 

{\rm[ii]} If $(\gA,\Phi,\shift{})$ is topologically weakly mixing,  then 
$\Specae=\{\bone\}$. 

{\rm[iii]} If $(\gA,\Phi,\shift{})$ is not topologically weakly mixing, and
 there is a $(\Phi,\shift{})$-ergodic measure with full support on $\gA$,
then $\Specae$ is nontrivial.
}
\bthmprf
 {[i]}  The proofs of {\bf(a,c,d,e)} are exactly as in 
 Lemma \ref{spectral.lemma}.
To see {\bf(b)}, suppose $f_1,f_2\in\Eigenspace$.  Then $g:=f_1/f_2$
is \ae-well-defined (because $|f_2|\ \eeequals{ae} \  1$)
 and $g\in \Eigenspace[\bone]$, where $\bone:=(1,...,1)$ [by part {\bf(a)}].
But if $(\gA,\Phi,\shift{})$ is transitive, then 
any element of $\Eigenspace[\bone]$ is \ae-constant by \cite[Thm 2.2]{KeynesRobertson}.  So $g \ \eeequals{ae} \ c $ for some constant $c\in\Cplx$; hence 
$f_1 \ \eeequals{ae} \ c\cdot f_2$.

[ii] is \cite[Thm 2.3]{KeynesRobertson}, and
[iii] is \cite[Thm 2.5 and Prop 2.6]{KeynesRobertson}.
\ethmprf

  An {\dfn eigenset} for $\gA$ is a collection $\set{f_\blam}{\blam\in\Specae}$
containing exactly one eigenfunction $f_\blam\in\Cae(\gA)$ for each $\blam\in\Specae$.
  
\Lemma{\label{eigenset}}
{  
\bthmlist
  \item $\Specae$ is countable {\rm(hence any eigenset is countable)}.

   \item Let $\sF=\{f_\blam\}_{\blam\in\Specae}$ be an eigenset for $\gA$.
There is a comeager, $(\Phi,\shift{})$-invariant subset
$\undgA=\undgA(\sF)\subseteq\gA$ such that, for any
$\undba\in\undgA$, any  $\blam\in\Specae$, 
any $\fz\in\ZD$, and $n\in\Natur$, \ 
\[
f_\blam\circ\Phi^n\circ\shift{\fz}(\undba)
\quad=\quad\lam_0^n\blam^{\fz} f_\blam(\undba). 
\]
  \item  If $\sF\subset\sC(\gA)$ {\rm [e.g. if $\Specae=\Spec$]}
then $\undgA=\gA$.
\ethmlist
}
\bthmprf
{\bf(a)} \ 
Any eigenset $\{f_\blam\}_{\blam\in\Specae}$ defines an \ae-continuous factor mapping
from $(\gA,\shift{})$
into a $\ZD[D+1]$-system $(\bT,\rho)$, where $\bT$ is a compact abelian group
and $\rho$ is a $\ZD[D+1]$-action by rotations of $\bT$.
We have $\Specae=\Spec[\bT,\rho]$.  Choosing one eigenfunction
for each $\blam\in\Spec[\bT,\rho]$, we get
an orthogonal basis of $\bL^2(\bT,\mu)$ (where $\mu$ is the Haar
measure on $\bT$).  But $\bL^2(\bT,\mu)$ is separable, so 
$\Spec[\bT,\rho]$ is countable.

{\bf(b)} \  Fix $\blam\in\Specae$.
For each $n\in\Natur$ and $\fz\in\ZD$, there is a comeager 
set $\gE_{(n;\fz)}\subseteq\gA$ such
 that, for any $\be\in\gE_{(n;\fz)}$, \  $f_\blam\circ\Phi^n\circ\shift{\fz}(\be)
=\lam_0^n\blam^{\fz} f_\blam(\be)$. 
Let $\gC_\blam:=\D\Intsct_{(n;\fz)\in\Natur\x\ZD} \gE_{(n;\fz)}$.
Then $\gC_\blam$ is comeager in $\gA$, and for every $n\in\Natur$,
$\fz\in\ZD$, and $\bc\in\gC_\blam$,
we have  $f_\blam\circ\Phi^n\circ\shift{\fz}(\bc)
=\lam_0^n\blam^{\fz} f_\blam(\bc)$.  
Let $\gB_\blam:=\D\Intsct_{(n;\fz)\in\Natur\x\ZD} 
\Phi^{-n}\shift{-\fz}(\gC_\blam)$; then $\gB_\blam$ is comeager,
and also 
$\Phi^n(\gB_\blam)\subseteq\gB_\blam$ and $\shift{\fz}(\gB_\blam)=\gB_\blam$,
for any $n\in\Natur$ and $\fz\in\ZD$.

  Finally, let $\undgA=\D\Intsct_{\blam\in\Specae} \gB_\blam$; then $\undgA$ is a
countable intersection [by {\bf(a)}] of comeager sets, and thus also
comeager.  Also, $\undgA$ is $(\Phi,\shift{})$-invariant because each
$\gB_\blam$ is $(\Phi,\shift{})$-invariant.

{\bf(c)}  If $f_\blam\in\sC(\gA)$, then $\gB_\blam=\gA$, because
$\gC_\blam=\gA$, because
$\gE_{(n;\fz)}=\gA$ for
all $n\in\Natur$ and $\fz\in\ZD$.  Thus, $\undgA=\gA$.
\ethmprf

\EXAMPLE{\label{X:eigenset}
(a) Let $\gS$  be as in Example \ref{X:defect.energy}(c).
Let $\sF:=\{f_1,f_{-1}\}$, where
$f_1:=\bone$ and where
 $f_{-1}:=f\in\Eigenspace[-1](\gS,\shift{})$ is from Example
\ref{X:ae.continuous}(a).  Then  $\undgS(\sF)=\gS\setminus\{\bar0\}$.

(b) If $\gT$ is as in Example \ref{X:ae.continuous}(b),
then $\undgT=\gT$ [because $\Specae[\gT]=\Spec[\gT]$].
}
 For any  $\dY\subset\ZD$ and $r\in\Natur$, let
$\dY(r):=\set{\fy\in\dY}{\dB(\fy,r)\subset\dY}$.
Recall that $\dY$ is {\dfn spacious} if $\dY(r)\neq\emptyset$ for all
$r\in\Natur$. A subshift $\gA\subset\AZD$
is {\dfn projectively transitive} if, for any
spacious $\dY\subset\ZD$ and  open $\gO\subset\gA$, the set
$\D\Union_{\fy\in\dY} \shift{-\fy}(\gO)$ is dense in $\gA$.

\EXAMPLE{
(a) If $\gA\subset\AZ$, then $\gA$
is projectively transitive iff both $\shift{1}$ and $\shift{-1}$ are
forward-transitive on $\gA$. 
[To see this, use Example \ref{X:projective.spacious}(a).]

(b) Suppose $\gA\subset\AZD$ has a $\shift{}$-ergodic measure $\mu$
with full support (i.e. $\mu[\gO]>0$ for any open $\gO\subseteq\gA$).
Then $\gA$ is projectively transitive.
To see this, let $\dY\subseteq\ZD$ be spacious and
let $\gO\subseteq\gA$ be open. 

\Claim{There is a subset $\gX\subseteq\D\Union_{\fy\in\dY}
\shift{-\fy}(\gO)$ with $\mu[\gX]=1$.}
\bclaimprf 
  For each $r\in\Natur$,
find $\fy_r\in\dY$ with $\dF_r:=\dB(\fy_r,r)\subset\dY$,
and define the function
\[
\alpha_r \quad:=\quad
\frac{1}{\#(\dF_r)} \sum_{\ff\in\dF_r} \chr{\gO}\circ\shift{\ff}.
\]
The sequence of sets $\{\dF_r\}_{r=1}^\oo$ is a \Folner sequence,
so the generalized Mean Ergodic Theorem 
\cite{OrnsteinWeissErgodic,Tempelman} says that
the sequence $\{\alpha_r\}_{r=1}^\oo$
converges to the constant function $\mu[\gO]$ in $\bL^2(\gA,\mu)$.
Thus, there is a subsequence $\{r_n\}_{n=1}^\oo$ and a set
$\gX\subseteq\gA$ with $\mu[\gX]=1$ such that $\D\lim_{n\goto\oo}
\alpha_{r_n}(\bx) = \mu[\gO]$ for all $\bx\in\gX$ \cite[Corol.2.32]{Folland}.
But $\gO$ is open, so $\mu[\gO]>0$ (because $\mu$ has full support on $\gA$). 
Thus, for any $\bx\in\gX$, there exists
$n\in\Natur$ such that $\alpha_{r_n}(\bx)>0$ (indeed, infinitely many such
$n\in\Natur$), which means there exists $\ff\in\dF_{r_n}$ such that
$\chr{\gO}\circ\shift{\ff}(\bx)=1$ which means $\shift{\ff}(\bx)\in\gO$.
But $\ff\in\dF_{r_n}\subset\dY$, so this means $\D\bx\in\Union_{\fy\in\dY}
\shift{-\fy}(\gO)$.\eclaimprf
Finally, $\gX$ is dense in $\gA$, because
$\mu[\gX]=1$ and $\mu$ has full support on $\gA$.
}
  Heuristically speaking, an eigenfunction $f$ detects some
underlying `rigidity' in the structure of $\gA$.  Thus, we don't need
to know every coordinate of $\ba\in\gA$ to evaluate $f(\ba)$; it
suffices to have information about some `large enough fragment' of
$\ba$.  This is the idea of the next lemma, 
where `large enough' means `spacious':

\Lemma{\label{ae.uniform.continuity}}
{  Let $\gA\subset\AZD$ be a projectively transitive subshift.
Let  $\sF:=\{f_\blam\}_{\blam\in\Specae}$ be an eigenset
and let $\undgA:=\undgA(\sF)$.
Let $\dY\subset\ZD$ be spacious.
\bthmlist
\item  
 There is a $(\Phi,\shift{})$-invariant,
comeager subset $\Good = \Good(\dY)\subseteq\undgA$ with the following
{\rm Extension Property:}

For any $\ba\in\Good$ and $\undba\in\undgA$,
 if $\undba_\dY=\ba_\dY$, then for every $\blam\in\Specae$, 
$f_\blam(\undba)=f_\blam(\ba)$.

\item If $\sF\subset\sC(\gA)$ {\rm [e.g. if $\Specae=\Spec$]}
then $\Good=\undgA=\gA$.

\item If $\dY'\subset\dY$, then $\Good(\dY')\subseteq\Good(\dY)$.

\item If $\dY':=\dY(R)$ for some $R\in\Natur$, then 
 $\Good(\dY')=\Good(\dY)$.

\item If $\Phi$ has radius $R$, and $\dY':=\dY(R)$, then
$\Phi[\Good(\dY)] \subseteq\Good(\dY')$.
\ethmlist
}
\bthmprf {\bf(a)} \ Fix $\blam\in\Specae$.

\setcounter{claimcount}{0}
\Claim{\label{ae.uniform.continuity.C1}
For any $\eps>0$, there exists  $r_\eps>0$ and a $(\Phi,\shift{})$-invariant,
comeager subset $\gG_\eps\subseteq\gA$ with the following property:
for any $\bg\in\gG_\eps$, there is some $\fy=\fy(\bg)\in\dY$
with $\dB(\fy,r_\eps)\subset\dY$ such that:
\beqn
\label{foobar}
\mbox{For any $\ba\in\gA$,}\qquad 
\statement{$\ba_{\dB(\fy,r_\eps)}=\bg_{\dB(\fy,r_\eps)}$}
\IMPLIES
\statement{$f_\blam\lb(\shift{\fy}[\ba]\rb) \closeto{\eps} f_\blam\lb(\shift{\fy}[\bg]\rb)$}\eeqn
}
\bclaimprf
 For all $r>0$, let \[
\gW_r(\eps):=\lb\{\bw\in\gA\ ; \ \forall \  \ba\in\gA, \ 
\statement{$\ba_{\dB(r)}=\bw_{\dB(r)}$}
\IMPLIES \rb.\lb.\statement{$f_\blam(\ba) \closeto{\eps} f_\blam(\bw)$}\rb\}.
\]
\subclaim{If $r$ is large enough, then $\gW_r(\eps)$ has nonempty interior.}
\bsubclaimprf
Let $\bw\in\gW_r(\eps/2)$, and let 
$\gC:=\set{\bc\in\gA}{\bc_{\dB(r)}=\bw_{\dB(r)}}$ be a cylinder neighbourhood
around $\bw$;
we will show that $\gC\subseteq\gW_r(\eps)$.
For any $\bc,\bc'\in \gC$, we have 
$f_\blam(\bc)\closeto{\eps}f_\blam(\bc')$, because
$f_\blam(\bc)\closeto{\eps/2}f_\blam(\bw) \closeto{\eps/2}
 f_\blam(\bc')$, because $\bc_{\dB(r)}=\bw_{\dB(r)}=\bc'_{\dB(r)}$.
Thus, $\bc\in\gW_r(\eps)$ for all $\bc\in\gC$, so $\gC\subseteq\gW_r(\eps)$.

  It remains to show that,  
if $r$ is large enough, then $\gW_r(\eps/2)$ is nonempty.  To see this, recall
that $f_\blam\in\Cae$, so $f_\blam$ has continuity points.
If $\ba\in\gA$ is any such continuity  point, then
$\ba\in\gW_r(\eps/2)$ if $r$ is large enough.
\esubclaimprf

Let $r_\eps:=r$; then $\dY(r)$ is nonempty and is also spacious.
Let $\gU_\eps$ be the (nonempty) interior of $\gW_r(\eps)$.
Let $\D \gO_\eps:=\Union_{\fy\in\dY(r)} \shift{-\fy}(\gU_\eps)$.
For any
$\bg\in\gO_\eps$, there exists $\fy\in\dY(r)$ with $\shift{\fy}(\bg)\in\gU_\eps
\subseteq\gW_r(\eps)$.
Then $\bg$ and $\fy$ satisfy eqn.(\ref{foobar}).

 Let $\gG_\eps := 
\D \Intsct_{\fz\in\ZD} \Intsct_{n=0}^\oo  \shift{-\fz}\Phi^{-n} (\gO_\eps)$.
Then $\gG_\eps$ is $(\Phi,\shift{})$-invariant.
Also, $\gG_\eps\subseteq\gO_\eps$, so for every $\bg\in\gG_\eps$ there
is some $\fy\in\dY$ satisfying eqn.(\ref{foobar}).
To see that $\gG_\eps$ is comeager,
let $\gU^{n;\fz}_\eps:=\shift{-\fz}\Phi^{-n}(\gU_\eps)$ for each $n\in\Natur$
and $\fz\in\ZD$.
Then $\gU^{n;\fz}_\eps$ is open because $\Phi^n$ and $\shift{\fz}$ are continuous and $\gU$ is open.
Thus, 
\[
\shift{-\fz}\Phi^{-n}(\gO_\eps) = \Union_{\fy\in\dY(r)} \shift{-\fz} \Phi^{-n}
 \shift{-\fy}(\gU_\eps) 
= \Union_{\fy\in\dY(r)}  \shift{-\fy}\shift{-\fz}\Phi^{-n}(\gU_\eps) 
= \Union_{\fy\in\dY(r)}  \shift{-\fy}(\gU^{n;\fz}_\eps)
\] 
is
open and dense in $\gA$ (because $\gA$ is projectively transitive). 
Thus, $\gG_\eps$ is a countable intersection of dense open sets, hence
 $\gG_\eps$ is comeager.
\eclaimprf
  Now, let $\D \gD_\blam:=\undgA \, \intsct \, \Intsct_{n=1}^\oo
\gG_{1/n}$.  Then $\gD_\blam$ is  $(\Phi,\shift{})$-invariant and 
comeager in $\gA$.

\Claim{\label{ae.uniform.continuity.C2}
If $\bd\in\gD_\blam$, \  $\undba\in\undgA$ and $\bd_{\dY}=\undba_{\dY}$, then
\ $f_\blam(\bd)=f_\blam(\undba)$.}
\bclaimprf
Fix $\eps>0$.  First, we claim that $f(\bd) \closeto{\eps} f(\undba)$.
To see this, find $n\in\Natur$ with  $1/n<\eps$.  Let $r:=r_{1/n}$ and
$\fy:=\fy(\bd)$ be as in Claim 1.  Then
\beqn
\label{ae.uniform.continuity.e1}
\blam^\fy f(\bd)
\quad\eeequals{(*)} \quad f\lb(\shift{\fy}[\bd]\rb) 
\quad  \closeto{\eps} \quad f\lb(\shift{\fy}[\undba]\rb)
\quad\eeequals{(\dagger)} \quad\blam^\fy f(\undba).
\eeqn
$(*)$ is because $\bd\in\undgA$ and
$(\dagger)$ is because $\undba\in\undgA$.
Finally, ``$\closeto{\eps}$'' is by Claim 1, because
$\bd\in\gD_\blam\subseteq\gG_{1/n}$, and because
$\bd_{\dB(\fy,r)}=\undba_{\dB(\fy,r)}$ because
$\bd_{\dY}=\undba_{\dY}$.  
But, $|\blam^\fy|=1$, so eqn.(\ref{ae.uniform.continuity.e1})
implies  $f(\bd) \closeto{\eps} f(\undba)$.

This argument works for any $\eps>0$.  Thus,
$f(\bd)=f(\undba)$.  
\eclaimprf
Let $\gC_0:= \D\Intsct_{\blam\in\Specae} \gD_\blam$.  Then
$\gC_0$ has the Extension Property, and 
is  $(\Phi,\shift{})$-invariant.
Finally, $\gC_0$ is
comeager, because Lemma \ref{eigenset}(a) implies that
$\gC_0$ is a countable intersection of comeager sets.  

 Thus, $\gC_0$ satisfies all the requirements of {\bf(a)}.  However, in
preparation for the proofs of {\bf(d)} and {\bf(e)} below, we must
refine this construction somewhat.
For each $r\in\Natur$, repeat the above construction to
obtain a $(\Phi,\shift{})$-invariant comeager set $\gC_r$ satisfying 
the Extension Property for $\dY(r)$.  Let $\Good:=\Intsct_{r=0}^\oo \gC_r$;
then $\Good$ is also comeager and  $(\Phi,\shift{})$-invariant, and
$\Good$ satisfies the Extension Property for $\dY$.

{\bf(b)} \ Repeat the construction in {\bf(a)} for $f_\blam\in\sC(\gA)$.
  In Claim 1, we have $\gO_\eps=\gA$ [because
if $f$ is continuous, then $f$ is uniformly continuous (because $\gA$
is compact), so that $\gU_\eps=\gW_r(\eps)=\gA$ if $r$ is large enough].
Thus, $\gG_\eps=\gA$.  Also, $\undgA=\gA$ by Lemma \ref{eigenset}(c).
Thus, $\gD_\blam=\gA$ for each $\blam\in\Spec$.  Thus, 
$\D\gC_0:=\Intsct_{\blam\in\Spec} \gD_\blam= \gA$.  Likewise,
$\gC_r=\gA$ for all $r\in\Natur$; hence $\Good=\gA$.

{\bf(c)} \  Suppose  $\dY'\subset\dY$, and
repeat the construction in {\bf(a)}
for $\dY'$. In Claim 1, $\gO_\eps'\subseteq\gO_\eps$, so
$\gG_\eps'\subseteq\gG_\eps$.  Thus, $\gD_\blam'\subseteq\gD_\blam$
for each $\blam$;  thus $\gC_0'\subseteq\gC_0$.  Likewise,
for each $r\in\Natur$, we have $\gC_r'\subseteq\gC_r$  because
$\dY'(r)\subseteq\dY(r)$.  Thus, $\Good(\dY')\subseteq\Good(\dY)$.

{\bf(d)}  If $\dY':=\dY(R)$ for some $R\in\Natur$, then for all
$r\in\Natur$, \ $\dY'(r)=\dY(r+R)$ [this is because $\dB(r+R)=\dB(r)+\dB(R)$].
Thus, for every $r\in\Natur$, \ $\gC_r'=\gC_{r+R}$.
Thus, 
\[
\D\Good(\dY') \quad = \quad \Intsct_{r=1}^\oo \gC'_r
\quad = \quad  \Intsct_{r=1}^\oo \gC_{r+R}
\quad \supseteq \quad  \Intsct_{r=1}^\oo \gC_{r}
\quad = \quad \Good(\dY) \quad \suuupseteq{(*)} \quad \Good(\dY'),\]
 where $(*)$ is by {\bf(c)}. Hence $\Good(\dY')=\Good(\dY)$.

{\bf(e)} \  We have $\Phi[\Good(\dY)] \ \subseteq \ \Good(\dY) \ 
\eeequals{(*)} \ \Good(\dY(R))$,
where $(*)$ is by {\bf(d)}
\ethmprf

\ignore{
{\bf(c)} \ If $\bc,\bc'\in\gC^n$ 
then there exist $\bb,\bb'\in\gC$ such that $\Phi^n(\bb)=\bc$ and 
$\Phi^n(\bb')=\bc'$.  If $\bc_{\dY^n}=\bc'_{\dY^n}$, then we can
choose $\bb$ and $\bb'$ such that $\bb_{\dY}=\bb'_{\dY}$.
Thus, for any $\blam\in\Specae$, if $f=f_\blam$, then
\[
f(\bc) \quad  = \quad 
f\circ \Phi^n(\bb) \quad \eeequals{(*)} \quad \lam_0^n f(\bb)
\quad \eeequals{(\dagger)} \quad \lam_0^n f(\bb') 
\quad \eeequals{(\diamond)} \quad
 f\circ\Phi^n(\bb') \quad = \quad f(\bc').
\]
$(*)$ is because $\bb\in\undgA$ and $(\diamond)$ is 
because $\bb\in\undgA$. 
Finally, $(\dagger)$ is because $\bb,\bb'\in\gC$ and  $\bb_{\dY}=\bb'_{\dY}$.}

\EXAMPLE{\label{X:ae.uniform.continuity}
Let $\sA=\{0,1\}$.
Let $\dY\subset\Zahl$ be spacious [see Example \ref{X:projective.spacious}(a)].

(a) If $\gS\subset\AZ$  and $\sF=\{f_1,f_{-1}\}$
are as in Example \ref{X:eigenset}(a),
then $\h{\gS}(\dY)=\gS\setminus\{\bar0\}$.

(b) If $\gT\subset\AZ$ is as in Example \ref{X:ae.continuous}(b), then $\h{\gT}(\dY)=\gT$, by Lemma \ref{ae.uniform.continuity}(b).
}

\Remark{
The constructions of $\undgA(\sF)$ in Lemma \ref{eigenset}(a) and
$\Good(\dY)$ in Lemma \ref{ae.uniform.continuity}(a)
depend upon the eigenset $\sF$,
because any \ae-equivalence class in $\Eigenspace$
 could contain uncountably many functions, each pair of which differ
on a meager subset of $\gA$.  Two eigenfunctions
$f_\blam,f'_\blam\in\Eigenspace$ could thus yield two sets 
$\Good(\dY)$ and $\Good'(\dY)$
whose symmetric difference was meager.
}

  Fix an eigenset
$\sF$, and let $\undgA:=\undgA(\sF)$ be as in Lemma \ref{eigenset}(a). Define
\newcommand{\tlundgA}{\widetilde{\undgA}}
\[
 \tlundgA:=\set{\ba\in\tlgA}{\mbox{$\ba_\dX \in \undgA_\dX$, for
every $r>0$ and projective component $\dX$ of $\unflawed(\ba)$}}.
\]
  Heuristically, `almost all' elements of $\tlgA$ are in $\tlundgA$
(because almost all elements of $\gA$ are in $\undgA$).
For any spacious $\dY\subseteq\ZD$, we define $\Good[\dY]:=\set{\ba_\dY}{\ba\in
\Good(\dY)}\subset\sA^{\dY}$, where $\Good(\dY)$ is as in 
Lemma \ref{ae.uniform.continuity}.
For any collection  $\dY_1,\ldots,\dY_N$ of disjoint  spacious subsets
of $\ZD$, we define
\[
\hgA(\dY_1,\ldots,\dY_N) \ \ := \ \ 
\set{\ba\in\tlundgA}
{\forall \ n\in\CC{1...N}, \ \ba_{\dY_n}\in\Good[\dY_n]}.
\]
If $\ba\in\hgA(\dY_1,\ldots,\dY_N)$, and $\blam\in\Specae$, then for
 each $n\in\CC{1...N}$, \ Lemma \ref{ae.uniform.continuity}(a) says
 there is a unique value $c_n\in\Cplx$ such that $f_\blam(\undba)=c_n$
 for any `extension' $\undba\in\undgA$ with
 $\undba_{\dY_n}=\ba_{\dY_n}$.  We can thus define
 $f_\blam(\ba_{\dY_n}):=c_n$.  For any $n,m\in\CC{1...N}$, we then
 define $\del_{n,m}(\blam):=
 f_\blam(\ba_{\dY_n})/f_\blam(\ba_{\dY_m})$.  Note that we make no
 assumption here about the relationship between the sets
 $\dY_1,\ldots,\dY_N$ and the projective components of
 $\unflawed(\ba)$.  Clearly, if $\ba\in\gA$, then
 $\del_{n,m}(\blam)=1$ for all $n,m\in\CC{1..N}$. However, if
 $\del_{n,m}(\blam)\neq 1$, then $\ba$ must be defective, and $\dY_n$
 and $\dY_m$ must lie in different projective components of
 $\unflawed(\ba)$, by part {\bf(d)} of the next result:

\Theorem{\label{ae.spec}}
{
Let $\Phi:\AZD\into\AZD$ be a CA and let $\gA\subset\AZD$ be a 
projectively transitive subshift
with $\Phi(\Sft)=\Sft$.  Fix an eigenset $\sF=\{f_\blam\}_{\blam\in\Specae}$.
Let $\dY_1,\ldots,\dY_N$  be disjoint  spacious subsets
of $\ZD$ and let $\ba\in\hgA(\dY_1,\ldots,\dY_N)$.  Then
\bthmlist
  \item $\del_{n,m}(\blam)$ is well-defined independent of the choice
of $f_\blam\in\Eigenspace$.
  \item $\del_{n,m}:\Specae\into\dT$ is a group homomorphism;
i.e. $\del_{n,m}\in\hSpecae$.
  \item  $\del_{n \ell}(\blam) = \del_{nm}(\blam)\del_{m \ell}(\blam)$ for any  $n,m,\ell\in\CC{1...N}$.
  \item If $\del_{nm}\neq \boldsymbol{1}$, then $\ba$ has a projective
domain boundary between $\dY_n$ and $\dY_m$.

  \item If $\blam\in\RatSpec$, then this $\del_{n,m}(\blam)$ is the
same as the one in {\rm Theorem \ref{gap.prop}.} 
\ethmlist
}
\bthmprf
 {\bf(a)} is by Lemma \ref{keynes.robertson}[i](b).  Part {\bf(b)} is by
Lemma \ref{keynes.robertson}[i](a), as in Claim 1[c] of Theorem \ref{gap.prop}.
Part {\bf(c)} is true by definition of $\del_{nm}$, \
and {\bf(e)} is straightforward.

For {\bf(d)}, suppose by contradiction that
$\dY_n$ and $\dY_m$ lie in the same projective
component $\dX$ of $\unflawed(\ba)$.
Find some extension $\undba\in\undgA$ such that
$\undba_\dX=\ba_\dX$ (this $\undba$ exists because $\ba\in\tlundgA$).
Note that $\undba$ is
also an extension of $\ba_{\dY_n}$ and of  $\ba_{\dY_m}$.
Thus, for every $\blam\in\Specae$, we have
$f_\blam(\ba_{\dY_n})=f_\blam(\undba)=f_\blam(\ba_{\dY_m})$
by definition, so $\del_{nm}(\blam)=f_\blam(\ba_{\dY_n})/f_\blam(\ba_{\dY_m})=1$.

So, if $\del_{nm}(\blam)\neq 1$ for any $\blam\in\Specae$,
 then $\dY_n$ and $\dY_m$ must be in
different  projective components of $\unflawed(\ba)$, which means that
$\ba$ has a projective domain boundary between $\dY_n$ and $\dY_m$.
\ethmprf

 If $\del_{n,m}(\blam)\neq 1$, then we say $\ba$ has a {\dfn projective
 $(\gA,\Phi)$-dislocation} between $\dY_n$ and $\dY_m$.  (When $\gA$
 and $\Phi$ are clear, we will just call this a ``projective
 dislocation'').  The matrix   $\Del_{\ba}  :=  [\del_{nm}]_{n,m=1}^N$
is called the {\dfn displacement matrix} of $\ba$ with respect to
the sets $(\dY_1,\ldots,\dY_N)$.

\EXAMPLE{\label{X:irrational.dislocation}
(a) Let $\gS\subset\AZ$ 
 and $\sF=\{f_1,f_{-1}\}$ be as in Examples \ref{X:eigenset}(a)
and \ref{X:ae.uniform.continuity}(a).
If $\bs\in\tl\gS$ is as in Example \ref{X:defect.energy}(c),
then $\bs$ has a dislocation, with displacement $\del_\ba=1\in\Zahlmod{2}$.

(b) Let $\sA=\{0,1\}$,  let $\lam:=(\sqrt{5}-1)/2$
and let $\gT\subset\AZ$ be the `Fibonacci' Sturmian subshift \cite[\S5.4.3]{Fogg} generated  by $\lam$ [see Example \ref{X:ae.continuous}(b)].
\[
\begin{array}{r}
\mbox{A typical point is} \ \ 
(\ldots010010010100101001001\ 0100100101001010010010100\ldots).\\
\mbox{Instead, let} \ \ba \ := \ 
(\ldots010010010100101001001.1001001010010100100101001\ldots).
\end{array}
\]
 Then $\ba$ has a dislocation at the decimal point, with displacement 
$\del_\ba=1\in\Zahl$.

(c) Let $\Dom$ be as in Example \ref{X:codim.one}(c).
Despite appearances, 
the domain boundary in Figure \ref{fig:codim}(D) is {\em not} a
dislocation, because $(\Dom,\shift{})$ is topologically mixing 
\cite[Lemma 2.1]{Ein},
so $\Specae[\Dom,\shift{}]$ is trivial by Lemma
\ref{keynes.robertson}[ii].  Instead, this is 
a `gap' defect; see \cite[Example 2.14(c)]{PivatoDefect2}.
}

Our last major result, analogous to
Theorem \ref{persistent.dislocation}, is that
projective dislocations are $\Phi$-persistent
defects:

\Theorem{\label{persistent.proj.disl}}
{
 Let $\Phi:\AZD\into\AZD$ be a CA of radius $r\geq 0$.
Let $\gA\subset\AZD$ be a projectively transitive
subshift with $\Phi(\gA)=\gA$.
 Let $\dY_1,\ldots,\dY_N\subset\ZD$ be disjoint spacious sets.
For all $n\in\CC{1...N}$, let $\dY'_n:=\dY_n(r)$.
Let $\ba\in\tlgA$ and $\ba':=\Phi(\ba)$.
\bthmlist
  \item  If $\ba\in\hgA(\dY_1,\ldots,\dY_N)$, then 
$\ba'\in\hgA(\dY'_1,\ldots,\dY'_N)$.

  \item  If $\ba$ has a projective dislocation, then so does $\ba'$.
Indeed, the $(\dY'_1,...,\dY'_N)$-displacement matrix of $\ba'$
is equal to the $(\dY_1,...,\dY_N)$-displacement matrix of $\ba$.
\ethmlist
}
\bthmprf  {\bf(a)} Follows immediately from 
Lemma \ref{ae.uniform.continuity}(e).

{\bf(b)}  Fix $n,m\in\CC{1..N}$.  Let $\undba_n\in\undgA$ and $\undba_m\in\undgA$
be extensions of $\ba_{\dY_n}$ and $\ba_{\dY_m}$, 
respectively. 
Let $\undba_n':=\Phi(\undba_n)$
and $\undba_m':=\Phi(\undba_m)$; then 
$\undba'_n$ and $\undba'_m$
are also in $\undgA$ 
[because $\Phi(\undgA)\subseteq\undgA$ by Lemma \ref{eigenset}(b)],
and are extensions of $\ba'_{\dY'_n}$ and $\ba'_{\dY'_m}$, 
respectively.  Thus, for any $\blam=(\lam_0;\lam_1,\ldots,\lam_D)\in\Specae$,
we have
\[
\del'_{nm}(\blam) \quad\eeequals{(\dagger)}\quad
\frac{f(\undba'_n)}{f(\undba'_{m})}
\quad\eeequals{(*)} \quad
\frac{\lam_0 f(\undba_n)}{\lam_0 f(\undba_m)}
\quad=\quad
\frac{f(\undba_n)}{f(\undba_m)}
\quad\eeequals{(\ddagger)}\quad
\del_{nm}(\blam).
\]
$(*)$ is because $f_\blam(\undba'_n) = f_\blam\circ\Phi(\undba_n)=\lam_0 f_\blam(\undba_n)$, because $\undba_n\in\undgA$.
Likewise $f_\blam(\undba'_m) = \lam_0 f_\blam(\undba_m)$
because $\undba_m\in\undgA$.
$(\ddagger)$ is by definition of $\del_{nm}(\blam)$,
and is independent of the choice of $\undba_n$ and $\undba_m$, by
Lemma \ref{ae.uniform.continuity}(a). Likewise
$(\dagger)$ is by definition of $\del'_{nm}(\blam)$.  

This holds for all $n,m\in\CC{1...N}$ and $\blam\in\Specae$.
\ethmprf

\EXAMPLE{
The essential dislocation in Example \ref{X:irrational.dislocation}(a)
is $\ECA{18}$-persistent. Figure \ref{fig:defect.intro}(a) shows
the long-term evolution of such defects.
}

\NRemark{\label{weak.mixing.remark}
 For one-dimensional SFTs, topological weak mixing implies mixing.
However, for other 
one-dimensional symbolic dynamical systems (e.g. rank one systems)
this is not true.  If $\gA$ is topologically weakly mixing, but
not  mixing, then $\gA$ admits no dislocations 
(by Lemma \ref{keynes.robertson}[ii]), but
must still admit other essential defects
[by Example \ref{X:removable.finite.defects}(a)], which are
undetectable by the spectral invariants developed here.
Are there other spectral invariants which  detect these
defects?}

\subsection*{Conclusion}

  We have used spectral theory to explain the  persistence
and interaction of domain boundaries in cellular automata.  However,
many questions remain.
\blist
\projdisl{\item
 For one-dimensional SFTs, topological weak mixing implies mixing.
However, for other 
one-dimensional symbolic dynamical systems (e.g. rank one systems)
this is not true.  If $\gA$ is topologically weakly mixing, but
not  mixing, then $\gA$ admits no dislocations 
(by Lemma \ref{keynes.robertson}[ii]), but
must still admit other essential defects
[by Example \ref{X:removable.finite.defects}(a)], which are
undetectable by the spectral invariants developed here.
Are there other spectral invariants which can detect these
defects?}

\item  Domain boundaries also emerge in coupled-map lattices
\cite{Kaneko88a,Kaneko88b,Kaneko89c}; \cite[\S8.2.4]{Ilachinski}.
Can analogous spectral invariants be developed in this context?

\item \label{CQ:blight}
 In most of our examples (e.g. ECAs \#54, \#62, \#110, and
\#184), the defects remain bounded in size, and act like `particles' 
\cite{PivatoDefect0}.  In general, however, defects may grow over
time like `blights' which invade the whole lattice.  What are
necessary/sufficient conditions for the defect to remain bounded?

\item \label{CQ:natural} Aside from the aforementioned ECAs, there are relatively few
known examples of `naturally occuring' defect dynamics in CA, and none
in $\ZD$ with $D\geq 2$.  It is easy to contrive
artificial examples, but this generally does not yield any surprises.
Are there nontrivial
examples of defect dynamics in multidimensional cellular automata?
Can we find them without blindly searching
the (vast) space of possible rules?

\item  \label{CQ:fixed.image}
 If $\gA\subset\AZD$, and there is a CA $\Phi$ and $n\in\Natur$ with
$\Phi^n(\AZD)\subseteq\gA\subseteq\Fix{\Phi}$, then $\gA$ admits no
essential defects.  The converse is also true, when $\gA$ is a
one-dimensional sofic shift with a $\shift{}$-fixed point
\cite{Maass}.  Is the converse true in higher dimensions?  

\item  \label{CQ:slow.converge} Even when $\gA$ admits essential defects, 
\Kurka and Maass \cite{KuMa00,KuMa02,Kur03,Kur05} have described how a
one-dimensional CA can `converge in measure' to $\gA$ through a gradual
process of defect coalescence/annihilation.  Given a subshift
$\gA\subset\AZD$, is it possible to build a CA which converges to
$\gA$ in this sense?

\item \label{CQ:condense}
The defect dynamics in ECAs \#18, \#54, \#62, \#110, and \#184 were
easy to discover by accident, because each CA contains a `condensing'
subshift $\gA$, such that generic initial conditions rapidly
`condense' into sequences containing relatively few defects separated
by long, $\gA$-admissible intervals.  What are necessary/sufficient
conditions for the existence of such a condensing subshift?  (This
rapid primordial condensation is not the same as the long-term
convergence in question \#\ref{CQ:slow.converge}, but the
two may be related.)  This question is closely related to question
\#\ref{CQ:blight}, because condensation should prevent defects from
growing.  Also, it relates to question \#\ref{CQ:natural}, because a
characterization of CA with condensing subshifts might yield
 nontrivial examples of defect dynamics.
\elist

Finally, we remark that the spectral invariants in this paper are only
applicable to defects of codimension one (i.e. domain boundaries).  In
a companion paper
\cite{PivatoDefect2}, we develop algebraic invariants for defects
of codimensions two or more.

%\ignore
{
\parag{Acknowledgements:}   This paper was mostly written during a research leave at Wesleyan University.  I am grateful to Ethan Coven, Adam Fieldsteel, Mike and Mieke Keane, and Carol Wood for their generosity and hospitality. 

{%\scriptsize
%\tiny\twocolumn
%\bibliographystyle{alpha}
\bibliographystyle{fundam}
\bibliography{bibliography}
}
}
\end{document}